\documentclass[ba,preprint]{imsart}

\RequirePackage{amsthm,amsmath,amsfonts,amssymb}
\RequirePackage[numbers,sort&compress]{natbib}
\RequirePackage[colorlinks,citecolor=blue,urlcolor=blue,backref=page]{hyperref}
\RequirePackage{graphicx}

\usepackage[utf8]{inputenc} 
\usepackage[T1]{fontenc}    
\usepackage{url}            
\usepackage{booktabs}       
\usepackage{amsfonts}       
\usepackage{amsmath}
\usepackage{amsthm}
\usepackage[algoruled,vlined]{algorithm2e}
\usepackage{mathtools}
\usepackage{nicefrac}       
\usepackage{microtype}      
\usepackage{fancyhdr}       
\usepackage[export]{adjustbox}
\usepackage{aligned-overset}
\graphicspath{{media/}}     
\usepackage{enumitem}

\usepackage{aliascnt}
\usepackage[capitalise,noabbrev]{cleveref}

\usepackage{xspace}
\usepackage{xcolor}

\allowdisplaybreaks    

\let\originalleft\left
\let\originalright\right
\renewcommand{\left}{\mathopen{}\mathclose\bgroup\originalleft}
\renewcommand{\right}{\aftergroup\egroup\originalright}

\pubyear{2026}
\arxiv{2605.26881}
\volume{TBA}
\issue{TBA}
\firstpage{1}
\lastpage{1}

\startlocaldefs
\theoremstyle{plain}

\newtheorem{theorem}{Theorem}[section]
\newaliascnt{lemma}{theorem}
\newtheorem{lemma}[lemma]{Lemma}
\aliascntresetthe{lemma}
\theoremstyle{definition}
\newaliascnt{definition}{theorem}

\aliascntresetthe{definition}

\theoremstyle{remark}

\theoremstyle{plain}
\newaliascnt{corollary}{theorem}
\newtheorem{corollary}[corollary]{Corollary}
\aliascntresetthe{corollary}
\newaliascnt{proposition}{theorem}
\newtheorem{proposition}[proposition]{Proposition}
\aliascntresetthe{proposition}

\theoremstyle{definition}

\newtheorem{innercustomthm}{Assumption}
\newenvironment{assumption}[1]
  {\renewcommand\theinnercustomthm{#1}\innercustomthm}
  {\endinnercustomthm}

\def\disclosurename{Disclosure of tool use}
\newenvironment{disclosure}[1][\disclosurename]{\begin{acks}[\disclosurename]}{\end{acks}}

\newcommand{\supptitle}{\title}
\endlocaldefs

\begin{document}

\begin{frontmatter}
\pdftitle{Robust Ensemble Kalman Filtering Under Observation Noise Misspecification}

\title{Robust Ensemble Kalman Filtering Under Observation Noise Misspecification}
\runtitle{Robust EnKF Under Observation Misspecification}

\begin{aug}
\author[A]{\fnms{Hans}~\snm{Reimann}\ead[label=e1]{reimann@math.uni-heidelberg.de}}
\and
\author[B]{\fnms{Sebastian}~\snm{Reich}\ead[label=e2]{sebastian.reich@uni-potsdam.de}}
\address[A]{Institute for Mathematics,
Heidelberg University\printead[presep={,\ }]{e1}}

\address[B]{Institute for Mathematics,
University of Potsdam\printead[presep={,\ }]{e2}}

\runauthor{H. Reimann and S. Reich}
\end{aug}

\begin{abstract}
We address observation noise misspecification in Bayesian filtering of dynamical systems via generalised Bayesian inference. Mismatch in tail decay between the true data generating process and an assumed observation model, often indicated via outliers, can strongly impact Bayesian updates in Kalman filters and their ensemble approximations. Adjusting information processing in the analysis step by employing diffusion score matching achieves robustness while maintaining closed-form Gaussian updates and good uncertainty quantification. We provide theoretical properties of this novel Kalman filter in linear Gaussian state space systems covering conjugacy of the analysis step, robustness, covariance stability, and tuning. We extend the approach to non-linear signal dynamics by adapting ensemble Kalman filter variants. Simulation studies on 2D target tracking and the Lorenz-63 and 40D Lorenz-96 systems highlight a critical trade-off between robustness and stability in filtering.
\end{abstract}



\end{frontmatter}

\section{Introduction}\label{sec: intro}

Ensemble Kalman filters (EnKFs) \cite{evensen1994sequential,evensen2009data} are standard practice in data assimilation for high-dimensional, non-linear systems such as in numerical weather prediction \cite{houtekamer2016review}. Via modifications such as covariance inflation and localisation, EnKFs can remain stable even with ensemble sizes much smaller than state dimension \cite{evensen2022data}. Since the analysis step employs Gaussian approximation, it inherits the sensitivity of the Kalman filter to observation noise misspecification, e.g. via departure from assumed Gaussian observation noise (see e.g. \cite[Theorem 3.2]{duran2024outlier}).

Robust state estimation of dynamical systems and robust Kalman filtering have been recurring topics in various applied fields since the seminal works \cite{yatawara1991kalman,xie1994robust}. We refer to \cite{boustati2020generalised,duran2024outlier} for overviews of robust particle filters and Kalman filters. The work at hand relates most closely to the weighted observation likelihood filter (WoLF) in \cite{duran2024outlier}.

We approach the problem of robust filtering by adjusting information processing via generalised Bayesian inference (GBI) as pioneered in \cite{bissiri2016general,jewson2018principles}. GBI was first explored for Bayesian filtering via the generalised particle filter in \cite{boustati2020generalised} and only recently generalised Kalman filter variants were proposed in \cite{duran2024outlier,reimann2024towards}. While Bayesian filtering and Gibbs posteriors (closely related to GBI; see \cite{mcallester1998some}) have been studied for much longer, their combination here relies on the introduction of diffusion score matching (DSM) in \cite{barp2019minimum} and its use for GBI in \cite{altamirano2023robust,altamirano2023robust2}. The established conjugacy for Gaussian prior--likelihood pairs therein motivated investigating DSM-based Kalman filter variants.

Generalised KF variants employing alternative learning measures can provide a way to balance robustness and information intake in a dynamic, a priori way at the cost of optimality in the well-specified case. In contrast to other methods based on detect-and-reject schemes or post hoc covariance inflation, the proposed diffusion score matching Kalman filter and its ensemble approximations redistribute information via observation based rescaling of the observation precision by a factor in $(0,2]$ and a shift in the innovation.

Other works, such as \cite{provost2023adaptive,tang2024generalized}, often make alternative explicit noise assumptions, e.g. via heavy-tailed, t-distributed observation noise. We instead assume Gaussian observation noise to be justified up to influential instances of deviation between the true data generating process and the assumed model. Specifying assumptions on this true data generating process is a central point of investigation and discussion.

Robust particle filters, such as the generalised particle filter in \cite{boustati2020generalised}, address non-linear dynamics and easily adjust to non-Gaussian noise, yet suffer from a curse of dimensionality and require a much larger number of particles than EnKFs (see e.g. \cite{van2019particle}). 
GBI-based Kalman filters such as the DSM KF maintain computational tractability of EnKFs via a closed-form Gaussian parameter update. Additionally, they can account for departure from Gaussian observation noise with theoretical guarantees on robustness and stability. GBI is a comparatively recent and growing topic. Developing Kalman filter variants alongside lets stochastic filtering benefit from new insights in this field.

\paragraph{Structure and contributions.}
After notation and setup, \cref{sec: DSM step} proves conjugacy of the DSM posterior in the linear Gaussian inverse problem and its global bias robustness for suitable weight kernels. \Cref{sec: DSM KF} formulates the DSM Kalman filter and proves weak stochastic boundedness of the analysis covariance given innovations with bounded second moments, bounded influence of the observation on the analysis mean, and, for a proposed default tuning, asymptotically unbiased analysis precision in the observation dimension limit under a well-specified observation marginal. \Cref{sec: ens nl} derives the stochastic DSM ensemble Kalman filter (EnKF), the deterministic DSM ensemble square root filter (ESRF), and the DSM local ensemble transform Kalman filter (LETKF) via standard arguments. \Cref{sec: simulation study} illustrates the theoretical properties in 2D target tracking and compares the regular, DSM, and WoLF  variants of the EnKF and LETKF on the Lorenz-63 and 40D Lorenz-96 models for varying degrees of misspecification. \Cref{sec: discussion} discusses the assumptions on the true data generating process and contextualises results via a trade-off between robustness and stability.

\section{Setup and motivation}

The Bayesian learning can be expressed via
\begin{equation}\label{eqn: KL Bayes}
    p(x_n|y_{1:n})= \frac{p(x_n|y_{1:(n-1)})\cdot p(y_n|x_n)}{p(y_n)}
    = p(x_n|y_{1:(n-1)})\cdot\exp\left[-\widehat{\mathrm{KL}}[p(y_n)\|p(y_n|x_n)]\right]
\end{equation}
with $p(y_n)$ short for the observation marginal $p(y_n|y_{1:(n-1)})$ and $\widehat{\mathrm{KL}}[p(y_n)\|p(y_n|x_n)]=-\log\left[\frac{p(y_n|x_n)}{p(y_n)}\right]$ the one-sample Monte Carlo estimator of the Kullback--Leibler divergence between the conditional and marginal distribution of the observation,
\begin{equation}\label{eqn: KL expectation}
    \mathrm{KL}[p(\,\cdot\,)\|p(\,\cdot\,|x_n)]=\mathbb{E}_{p}\left[-\log\left[\frac{p(Y_n|x_n)}{p(Y_n)}\right]\right]=-\int_{\mathbb{R}^{d_Y}}p(y_n)\log\left[\frac{p(y_n|x_n)}{p(y_n)}\right]\mathrm{d}y_n\,.
\end{equation}
Classical Bayesian learning can be understood in an optimal information processing sense (see e.g. \cite{zellner1988optimal,kullback1997information} for details).

A popular application of Bayesian learning is inference in stochastic filtering. We investigate the established linear Gaussian state space model (LGSS) as a popular special case. Notation and arguments follows \cite{reich2015probabilistic} with some adjustments. Multivariate Gaussian random variables are denoted in covariance form via $p(x)=n(x;m,P)$ with covariance matrix $P$ and mean vector $m$ and information form via $p(x)=n^{-1}(x;\theta,J)$ with precision matrix $J=P^{-1}$ and potential $\theta=Jm$. $\nabla_x\cdot f(x)=\langle\nabla,f\rangle$ denotes the divergence operator on $f$ and $\nabla_x f(x)$ denotes the gradient of $f$.

Let $(\Omega,\mathcal{F},\mathbb{P})$ be the underlying probability space.
\begin{assumption}{1.LGSS}[Linear Gaussian State Space System]\label{ass: LGSS} Let $n\in\mathbb{N}$, $d_X,d_Y\in\mathbb{N}$, and
\begin{align*}
    X_n&=A_nX_{n-1}+Q^{\frac{1}{2}}_nW_n\\
    Y_n&=H_nX_n+R^{\frac{1}{2}}_nV_n
\end{align*}
with 
\begin{itemize}
    \item $X_n:\Omega\to\mathbb{R}^{d_X}$ the $d_X$-dimensional signal random vector at time $n$,
    \item $Y_n:\Omega\to\mathbb{R}^{d_Y}$ the $d_Y$-dimensional observation random vector at time $n$,
    \item $W_n:\Omega\to\mathbb{R}^{d_X}$ and $V_n:\Omega\to\mathbb{R}^{d_Y}$ \emph{i.i.d.} standard Gaussian random vectors at time $n$ of corresponding dimensions,
    \item $A_n\in\mathbb{R}^{d_X\times d_X}$, $H_n\in\mathbb{R}^{d_Y\times d_X}$, and positive definite, symmetric $Q_n\in\mathbb{R}^{d_X\times d_X}$ and $R_n\in\mathbb{R}^{d_Y\times d_Y}$, and
    \item $p(x_0)=n(x_0;m^a_0,P^a_0)$, the initial Gaussian prior distribution at time $0$.
\end{itemize}
\end{assumption}
The celebrated Kalman filter solving the filtering problem in \cref{ass: LGSS} exploits that linear combinations of Gaussian random variables remain Gaussian plus conjugacy in the linear Bayesian inverse problem for Gaussian prior--likelihood pairs. Closed-form recursive update formulas of parameters then enable fast and accurate computation while maintaining full characterisation of the forecast, filtering, and marginal distributions (see \cite{reich2015probabilistic} for additional details).

We recall some key aspects of the regular KF. The forecast distribution is $p(x_n|y_{1:(n-1)})=n(x_n;m_n^f,P_n^f)$ with
\begin{itemize}
    \item forecast covariance $P_n^f=A_nP^a_{n-1}A_n^T+Q_n$ and
    \item forecast mean $m_n^f=A_nm^a_{n-1}$.
\end{itemize} 
We employ the forecast as prior distribution in the Bayesian inverse problem with the observation likelihood given by $p(y_n|x_n)=n(y_n;H_nx_n,R_n)$. Both are direct results of the system formulation in \cref{ass: LGSS}. The regular posterior of the linear Bayesian inverse problem is then given via 
\begin{equation*}
     p(x_n|y_{1:n}) \propto p(x_n|y_{1:(n-1)})\cdot p(y_n|x_n)
 \end{equation*}
with $p(x_n|y_{1:n})=n(x_n;m^a_n,P^a_n)$ and parameter update via
\begin{itemize}
    \item Kalman gain $K_n=P_n^fH_n^T\left[R_n+H_nP_n^fH_n^T\right]^{-1}$,
    \item analysis covariance $P^a_n=P_n^f-K_nH_nP_n^f$, and
    \item analysis mean $m^a_n=m_n^f-K_n\left(H_nm_n^f-y_n\right)$.
\end{itemize}

Classical Bayesian learning and the regular Kalman filter can be understood in an optimal information processing sense (see e.g. \cite{zellner1988optimal,kullback1997information} for details). Taking the perspective in \cref{eqn: KL Bayes}of the Bayesian update as statistical learning with respect to an information criterion, we investigate changing this discrepancy measure in generalised Bayesian filtering to obtain robustness under misspecification, i.e. when the true data generating process (DGP) $\pi(y_n)$ of $Y_n$ given $y_{1:(n-1)}$ differs from $p(y_n)$ (see e.g. \cite{knoblauch2018spatio, altamirano2023robust, altamirano2023robust2}).  

\section{The diffusion score matching analysis step}\label{sec: DSM step}

We replace the KL divergence in Bayes' formula as in \cref{eqn: KL expectation} with diffusion score matching \cite{barp2019minimum} as an estimator of diffusion Fisher divergence \cite{gong2021interpreting}. We extend preliminary results in \cite{reimann2024towards}.

Let $\widehat{\mathrm{DSM}}$ denote the one-sample diffusion score matching estimator replacing the one-sample KL estimator in \cref{eqn: KL Bayes}. Our main object of interest is therefore the corresponding generalised posterior,
\begin{equation}\label{eqn: DSM Bayes}
        p_\mathrm{DSM}(x_n|y_{1:n})\propto p(x_n|y_{1:(n-1)})\cdot\exp\left[-\widehat{\mathrm{DSM}}[\pi(y_n)\|p(y_n|x_n)]\right]\,.
\end{equation}
We assume the modelled observation marginal $p(y_n)$ to be an inaccurate representation, e.g. with respect to tail behaviour and frequency of supposed outliers. Accordingly, we express via $\pi(\,\cdot\,)$ in \cref{eqn: DSM Bayes} the true DGP. Diffusion score matching as introduced in \cite{barp2019minimum} is given via
\begin{equation}\label{eqn: DSM E}
    \mathrm{DSM}[\pi(\,\cdot\,)\|p(\,\cdot\,|x)]=\mathbb{E}_{\pi}\left[\left\|w^T(Y)s_{p(\,\cdot\,|x)}(Y)\right\|^2_2+2\nabla_Y\cdot\left(w(Y)w^T(Y) s_{p(\,\cdot\,|x)}(Y)\right)\right]
\end{equation}
with 
\begin{itemize}
    \item score function $s_{p(\,\cdot\,|x)}(y)=\nabla_y\log\left[p(y|x)\right]$ and
    \item pointwise invertible matrix-valued function $w:\mathbb{R}^{d_Y}\to\mathbb{R}^{d_Y\times d_Y}$.
\end{itemize}
    
DSM provides posterior robustness similar to the sense of Huber (see e.g. \cite{huber2004robust}). For a suitable $w$ (\cref{ass: kernel}), the DSM posterior is globally bias-robust (\cref{thm: bias rob}). This is not the case for the regular KL-based posterior (see e.g. \cite[Theorem 3.2]{duran2024outlier}). This is complemented by numerical tractability. Looking at \cref{eqn: DSM E}, diffusion score matching only contains the true DGP $\pi$ via its expectation. Accordingly, its Monte Carlo estimator only requires knowledge of $\pi$ via samples. Furthermore, results in \cite{altamirano2023robust,altamirano2023robust2} show the posterior \cref{eqn: DSM Bayes} has conjugacy properties for exponential family distributions including Gaussian prior--likelihood pairs. This suggests a recursive parameter update akin to the regular KF for the LGSS setting in \cref{ass: LGSS} (see \cite{reimann2024towards}).

DSM can be derived from diffusion Fisher divergence via an integration by parts argument akin to regular score matching with respect to regular Fisher divergence (see e.g. \cite{hyvarinen2005estimation}). This requires smoothness and boundary conditions on the true DGP in \cref{ass: regularity}. The one-sample DSM posterior in \cref{eqn: DSM Bayes} is well-defined for any observation $y\in\mathbb{R}^{d_Y}$ irrespective of $\pi$.
\begin{assumption}{2.R}[Regularity of the true DGP]\label{ass: regularity}
    The true DGP $\pi$ satisfies \\$\mathbb{E}_\pi\left[\left\|s_{p(\,\cdot\,|x)}(y)-s_{\pi(\,\cdot\,)}(y)\right\|^2\right]<\infty$ and $[\pi ww^Ts_{p(\,\cdot\,|x)}],[\nabla\cdot(\pi ww^Ts_{p(\,\cdot\,|x)})]\in L^1(\mathbb{R}^{d_Y})$ for observation likelihood $p(\,\cdot\,|x)$ and pointwise invertible matrix-valued $w:\mathbb{R}^{d_Y}\to\mathbb{R}^{d_Y\times d_Y}$.
\end{assumption}

\subsection{Gaussian conjugacy}

For stating conjugacy and the analysis parameter update, we drop the time-varying notation in \cref{ass: LGSS} for the remainder of this section.
\begin{assumption}{3.LIP}[Linear Inverse Problem]\label{ass: LIP}
    Under \cref{ass: LGSS} without the time index, let $p(x)=n(x;m^f,P^f)$ with $x\in\mathbb{R}^{d_X}$, $p(y|x)=n(y;Hx,R)$ with $y\in\mathbb{R}^{d_Y}$ and $m^f,P^f,H,R$ corresponding components of appropriate dimensions.
\end{assumption}

\begin{proposition}[Gaussian conjugacy of the DSM posterior]\label{thm: Conjugacy}
    Given the linear inverse problem in \cref{ass: LIP}. Choosing $w(y)=k(y)R^{\frac{1}{2}}$ with $k:\mathbb{R}^{d_Y}\to(0,1]$ differentiable, we obtain for the DSM posterior as given in \cref{eqn: DSM Bayes} that
    \begin{equation*}
        p_{\mathrm{DSM}}(x|y)=n(x;m^a,P^a)
    \end{equation*}
    with
    \begin{itemize}
        \item rescaled observation covariance $N(y)=\frac{1}{2k^2(y)}R$,
        \item shifted observation $\tilde{y}=y-2N(y)\nabla_{y}k^2(y)$,
        \item adjusted Kalman gain $\widetilde{K}(y)=P^fH^T\left[N(y)+HP^fH^T\right]^{-1}$,
        \item analysis covariance $P^a=P^f-\widetilde{K}(y)HP^f$, and
        \item analysis mean $m^a=m^f-\widetilde{K}(y)\left[Hm^f-\tilde{y}\right]$.
    \end{itemize}
    for an arbitrary fixed observation $y\in\mathbb{R}^{d_Y}$.
\end{proposition}

\begin{proof}
For the score function of the observation likelihood, observe
\begin{equation*}
    s_{p(\,\cdot\,|x)}(y)=\nabla_y\log[p(y|x)]=\nabla_y[-\frac{1}{2}(y-Hx)^TR^{-1}(y-Hx)]=-R^{-1}(y-Hx)\,.
 \end{equation*}
The diffusion matrix is chosen as $w(y)=k(y)R^{\frac{1}{2}}$ for $k:\mathbb{R}^{d_Y}\to(0,1]$. The weight kernel $k$ does not depend on $x$.

We obtain for the one-sample DSM Monte Carlo estimator that
\begin{align*}
    \widehat{\mathrm{DSM}}[\pi(\,\cdot\,)\|p(\,\cdot\,|x)]&=\|w^T(y)s_{p(\,\cdot\,|x)}(y)\|^2_2+2\nabla_{y}\cdot[w(y)w^T(y)s_{p(\,\cdot\,|x)}(y)]\\
    &\overset{+C}{=}(Hx)^Tk^2(y)R^{-1}Hx-2(Hx)^Tk^2(y)R^{-1}y+2(Hx)^T\nabla_{y}k^2(y)\\
    &=\frac{1}{2}(Hx)^TN(y)^{-1}Hx-(Hx)^TN(y)^{-1}\tilde{y}
 \end{align*}
with $N(y)^{-1}=2k^2(y)R^{-1}$, $\tilde{y}=y-N(y)\nabla_{y}2k^2(y)$ and $\overset{+C}{=}$ indicating equality up to an additive constant independent of $x$.

Given the prior in information form, so $p(x)=n^{-1}(x;\theta^f,J^f)$ and \cref{eqn: DSM Bayes} simplified for \cref{ass: LIP}, then
\begin{align*}
    p_\mathrm{DSM}(x|y)&\propto p(x)\exp\left(-\widehat{\mathrm{DSM}}[\pi(\,\cdot\,)\|p(\,\cdot\,|x)]\right)\\
    &\propto\exp\left(-\frac{1}{2}x^TJ^fx+x^T\theta^f\right)\exp\left(-\frac{1}{2}(Hx)^TN(y)^{-1}Hx+(Hx)^TN(y)^{-1}\tilde{y}\right)\\
    &=\exp\left(-\frac{1}{2}x^TJ^ax+x^T\theta^a\right)
 \end{align*}
 with
\begin{equation*}
        J^a\coloneqq J^f+H^TN(y)^{-1}H\quad\text{and}\quad\theta^a\coloneqq\theta^f+H^TN(y)^{-1}\tilde{y}\,,
\end{equation*}
the recursive parameter update in information form and $p_\mathrm{DSM}(x|y)=n^{-1}(x;\theta^a,J^a)$.

It remains to convert from information form to covariance form using the \emph{Sherman--Morrison--Woodbury} matrix inversion formula. For the covariance matrix,
\begin{equation*}
    P^a=(J^a)^{-1}=[J^f+H^TN(y)^{-1}H]^{-1}=P^f-\widetilde{K}(y)HP^f
\end{equation*}
with adjusted Kalman gain matrix
\begin{equation*}
    \widetilde{K}(y)=P^fH^T[N(y)+HP^fH^T]^{-1}\,.
\end{equation*}
Based on this, observe for the mean vector
\begin{align*}
    m^a&=P^a\theta^a=P^a[\theta^f+H^TN(y)^{-1}\tilde{y}]=m^f-P^aH^TN(y)^{-1}[Hm^f-\tilde{y}]\quad\mathrm{or}\\
    &=m^f-\widetilde{K}(y)(Hm^f-\tilde{y})\,.
\end{align*}
See \cite[Chapter 5.1]{reich2015probabilistic} for additional details. 
\end{proof}

In this work, we consider an adjusted inverse multi-quadratic (IMQ) kernel
\begin{equation}\label{eqn: IMQ kernel}
    k(y)=\left(1+\frac{\|y-Hm^f\|_{\Sigma^{-1}}^2}{q^2}\right)^{-\frac{1}{2}}
\end{equation}
with threshold parameter $q^2>0$ and squared Mahalanobis distance $\|y-Hm^f\|_{\Sigma^{-1}}^2$ standardising via the marginal covariance matrix $\Sigma=HP^fH^T+R$ (see e.g. \cite{mahalanobis2018generalized}).

\subsection{Posterior robustness}\label{sec: bias rob}

The notion of robustness considered here follows \cite{altamirano2023robust,altamirano2023robust2}. It shares close ties to classical robustness in the sense of Huber \cite{huber2004robust}, except that we do not consider point estimators but posteriors as estimators of distributions. The established concept of the influence function therefore needs adapting to that of the posterior influence function (PIF) capturing the effect of contamination on the posterior. We refer to \cite{ghosh2016robust,matsubara2022robust} for additional details.

We model misspecification via a true DGP $\pi$ that is an $\varepsilon$-contaminated version of the assumed observation marginal $p(y)=n(y;Hm^f,\Sigma)$. The true DGP enters the DSM loss in \cref{eqn: DSM E} only as the measure of the expectation, which the one-sample DSM posterior in \cref{eqn: DSM Bayes} replaces by an observation $y$. For a probability measure $\tilde{\pi}$ on $\mathbb{R}^{d_Y}$, denote
\begin{equation*}
    p_\mathrm{DSM}(x|\tilde{\pi})\propto p(x)\exp\left[-\mathrm{DSM}[\tilde{\pi}(\,\cdot\,)\|p(\,\cdot\,|x)]\right]
\end{equation*}
with $\mathrm{DSM}$ as in \cref{eqn: DSM E}. For $\tilde{\pi}=\delta_{y}$, the empirical measure of a single instance $y$ of the true DGP $\pi$, this recovers \cref{eqn: DSM Bayes}. We consider the $\varepsilon$-contamination
\begin{equation}\label{eqn: cont dgp}
    \tilde{\pi}_{\varepsilon,y_0}=(1-\varepsilon)\tilde{\pi}+\varepsilon\delta_{y_0}\quad\text{with}\quad y_0\in\mathbb{R}^{d_Y}\quad\text{and}\quad\varepsilon\in[0,1]\,.
 \end{equation}
If $\tilde{\pi}$ has density $p(y)$, \cref{eqn: cont dgp} models the misspecified true DGP and is discussed in \cref{sec: discussion}. The posterior influence function is then
\begin{equation*}
    \mathrm{PIF}(y_0,x,\tilde{\pi})=\left.\frac{\mathrm{d}}{\mathrm{d}\varepsilon} p_\mathrm{DSM}(x|\tilde{\pi}_{\varepsilon,y_0})\right\vert_{\varepsilon=0}\,.
\end{equation*}
We consider $p_\mathrm{DSM}(x|\tilde{\pi})$ globally bias-robust if there is a finite double uniform bound on the PIF,
\begin{equation}\label{eqn: global bias rob}
    \sup_{x\in\mathbb{R}^{d_X},\ y_0\in\mathbb{R}^{d_Y}}\left|\mathrm{PIF}(y_0,x,\tilde{\pi})\right|<\infty\,.
\end{equation}

The weight kernel $k$ was not further specified in \cref{thm: Conjugacy}. It introduces the novel degree of freedom to enable posterior robustness.

\begin{assumption}{4.$k$}[Weight kernel properties]\label{ass: kernel}
The kernel $k:\mathbb{R}^{d_Y}\to(0,1]$ is differentiable with $\sup_{y}k(y)\|y\|_2<\infty$, $\sup_{y}\|\nabla_yk^2(y)\|_2<\infty$, and $\sup_{y}\left|\nabla_y\cdot[k^2(y)y]\right|<\infty$.
\end{assumption}

The choice of IMQ kernel in \cref{eqn: IMQ kernel} is motivated via empirical success in \cite{altamirano2023robust,altamirano2023robust2,duran2024outlier} in similar contexts. It satisfies \cref{ass: kernel}, since $k^2(y)\|y-Hm^f\|^2_{\Sigma^{-1}}\leq q^2$ and $\nabla_yk^2(y)=-\frac{2}{q^2}k^4(y)\Sigma^{-1}(y-Hm^f)$.

\begin{theorem}[Global bias robustness of the DSM posterior]\label{thm: bias rob}
        Given the linear inverse problem in \cref{ass: LIP} and a weight kernel $k$ satisfying \cref{ass: kernel}. For $w(y)=k(y)R^{\frac{1}{2}}$ as in \cref{thm: Conjugacy}, the DSM posterior $p_\mathrm{DSM}(x|\tilde{\pi})$ is globally bias-robust as in \cref{eqn: global bias rob} for any probability measure $\tilde{\pi}$. In particular, this holds under the misspecification in \cref{eqn: cont dgp} with $p(y)$ the density of $\tilde{\pi}$, and for the single analysis step in \cref{thm: Conjugacy}.
\end{theorem}

\begin{proof}
We utilise \cite[Proposition B.1]{altamirano2023robust} and show
\begin{enumerate}
    \item $\underset{y\in\mathbb{R}^{d_Y}}{\sup}\left|\mathrm{DSM}\left[\delta_y\|p(\,\cdot\,|x)\right]\right|\leq\gamma(x)$ for $\gamma:\mathbb{R}^{d_X}\to[0,\infty)$,
    \item $\underset{x\in\mathbb{R}^{d_X}}{\sup}\ p(x)\gamma(x)<\infty$, and $\int_{\mathbb{R}^{d_X}}p(x)\gamma(x)\mathrm{d}x<\infty$.
\end{enumerate}

Investigating condition 1 via a triangle inequality,
\begin{align}\label{eqn: split} 
    \left|\mathrm{DSM}\left[\delta_y\|p(\,\cdot\,|x)\right]\right|&\leq\left\|w^T(y)s_{p(\,\cdot\,|x)}(y)\right\|^2_2+2\left|\nabla_{y}\cdot\left[w(y)w^T(y)s_{p(\,\cdot\,|x)}(y)\right]\right|\nonumber\\
    &\eqqcolon(i)^2+2(ii)\,,
\end{align}
for $y\in\mathbb{R}^{d_Y}$, covering both the contamination $y_0$ and the support of $\tilde{\pi}$.

Utilising \cref{ass: kernel} for controlling the terms in \cref{eqn: split}, we observe 
\begin{align*}
    (i)&=\left\|w^T(y)s_{p(\,\cdot\,|x)}(y)\right\|_2=\left\|w^T(y)R^{-1}(y-Hx)\right\|_2\\
    &\leq\left\|w^T(y)R^{-1}y\right\|_2+\left\|w^T(y)R^{-1}Hx\right\|_2\\
    &\leq \tilde{c}^{(i)}_1+\tilde{c}^{(i)}_2\|Hx\|_2\eqqcolon\gamma_1(x)\\
    &\quad\implies\gamma_1^2(x)=c^{(i)}_1\|Hx\|_2^2+c^{(i)}_2\|Hx\|_2+c^{(i)}_3\,.
\end{align*}
Similarly for the second component, we observe
\begin{align*}
    (ii)&=\left|\nabla_{y}\cdot\left[w(y)w^T(y)s_{p(\,\cdot\,|x)}(y)\right]\right|=\left|\nabla_{y}\cdot\left[k^2(y)(y-Hx)\right]\right|\\
    &\leq\left\|\nabla_{y}k^2(y)\right\|_2\|Hx\|_2+\left|\nabla_{y}\cdot\left[k^2(y)y\right]\right|\leq c^{(ii)}_1\|Hx\|_2+c^{(ii)}_2\eqqcolon\gamma_2(x)\,.
\end{align*}
In both cases, we choose appropriate positive constants independent of $x$ and $y$ based on the system components in \cref{ass: LIP} and the choice of $k$ satisfying \cref{ass: kernel}. Taking everything together,
\begin{equation*}
    \sup_{y\in\mathbb{R}^{d_Y}}\left|\mathrm{DSM}\left[\delta_y\|p(\,\cdot\,|x)\right]\right|\leq\gamma_1^2(x)+2\gamma_2(x)=c_1\|Hx\|^2_2+c_2\|Hx\|_2+c_3\eqqcolon\gamma(x)
\end{equation*}
for positive constants $c_1,c_2,c_3$.

Condition 2 and the boundedness of $p(x)$ required in \cite[Proposition B.1]{altamirano2023robust} follow directly for the Gaussian prior $p(x)$ in \cref{ass: LIP}.
\end{proof}

\section{The DSM Kalman filter and properties}\label{sec: DSM KF}

We employ the conjugacy in \cref{thm: Conjugacy} under \cref{ass: LGSS} to obtain a discrete time recursive parameter update. As with the regular Kalman filter, we maintain Gaussian distributions after each forecast and analysis step.

\begin{corollary}[The diffusion score matching Kalman filter]\label{thm: conjF}
Given the LGSS system in \cref{ass: LGSS}, the diffusion score matching filtering equations can be evaluated via recursive parameter updates and the resulting distributions are Gaussian:
\begin{itemize}
    \item $p(x_n|y_{1:(n-1)})=n(x_n;m^f_n,P^f_n)$ (the forecast distribution),
    \item $p(y_n|y_{1:(n-1)})=n(y_n;H_nm^f_n,\Sigma_n)$ (the observation marginal distribution), and
    \item $p_\mathrm{DSM}(x_n|y_{1:n})=n(x_n;m^a_n,P^a_n)$ (the analysis distribution).
\end{itemize}
The parameters are given by the diffusion score matching Kalman filter forecast and analysis steps.
\begin{itemize}[leftmargin=*]
    \item The forecast step (unchanged) is $P^f_n=A_nP^a_{n-1}A_n^T+Q_n$ and $m^f_n=A_nm^a_{n-1}$.
    \item The analysis step is
    \begin{center}
    \begin{minipage}[t]{0.45\linewidth}
    \begin{itemize}[leftmargin=*]
        \item $k^2_n(y_n)=\left(1+\frac{\|y_n-H_nm_n^f\|_{\Sigma_n^{-1}}^2}{q^2}\right)^{-1}$ with threshold parameter $q^2>0$,
        \item $N_n(y_n)=\frac{1}{2k_n^2(y_n)}R_n$,
        \item $\tilde{y}_n=y_n-2N_n(y_n)\nabla_{y_n}k_n^2(y_n)$,
    \end{itemize}
    \end{minipage}\hfill
    \begin{minipage}[t]{0.53\linewidth}
    \begin{itemize}[leftmargin=*]
        \item $\Sigma_n=H_nP^f_nH_n^T+R_n$,
        \item $\widetilde{K}_n(y_n)=P_n^fH_n^T\left[N_n(y_n)+H_nP_n^fH_n^T\right]^{-1}$,
        \item $P_n^a=P_n^f-\widetilde{K}_n(y_n)H_nP_n^f$, and
        \item $m^a_n=m_n^f-\widetilde{K}_n(y_n)\left[H_nm_n^f-\tilde{y}_n\right]$.
    \end{itemize}
    \end{minipage}
    \end{center}
\end{itemize}
The recursion is started with $m^a_0$ and $P^a_0$ given in \cref{ass: LGSS}.
\end{corollary}

The proof is omitted. The forecast follows via linearity of Gaussian distributions and the analysis step is a direct repeated application of \cref{thm: Conjugacy}. The DSM KF maintains the computational complexity of the regular KF with bottleneck in inverting a matrix of dimension $d_Y\times d_Y$.

The construction is easily adjusted to block-diagonal $R$ by replacing $k$ in $w(y)=k(y)R^{\frac{1}{2}}$ with a diagonal matrix of individual weight kernels $k_b$ corresponding to the blocks $R_b$ recovering $N_b(y)=\frac{1}{2k^2_b(y)}R_b$ for the $b$-th block. The IMQ weight kernel in \cref{eqn: IMQ kernel} is then adjusted accordingly via
\begin{equation*}
    k_b(y)=\left(1+\frac{\|\left[y-Hm^f\right]_b\|_{\Sigma_b^{-1}}^2}{q_b^2}\right)^{-\frac{1}{2}}
\end{equation*}
with $[y]_b$ the entries of $y$, $\Sigma_b$ the diagonal block of $\Sigma$, and individual tuning $q_b^2>0$ corresponding to the $b$-th block of $R$. The shifted observation adjusts accordingly to the individual components. Blocks not at risk of misspecification may be chosen constant via  $k_b(y)\equiv\frac{1}{\sqrt{2}}$ recovering the regular KF update there. We expect results for the single analysis step to translate, but do not investigate this further here.

\subsection{Stability of the DSM Kalman filter}\label{sec: stability}

Complementary to robustness, we study the long-time behaviour of the DSM KF. The filter differs from the regular KF through the rescaled observation covariance $N_n(y_n)$ and the shifted observation $\tilde y_n$. For the IMQ kernel in \cref{eqn: IMQ kernel}, both depend on the filter state through $m^f_n$ and $\Sigma_n$. We therefore control the inflation of $N_n(y_n)$ via the innovation term. To account for the randomness introduced in the DSM analysis step, we consider stability in the sense of \cite{solo1996stability} with respect to weakly stochastically bounded components, i.e. $\lim_{\delta\to\infty}\sup_{n}\mathbb{P}\left(\|M_n\|>\delta\right)=0$ for a sequence of random matrices $(M_n)_{n\in\mathbb{N}}$.

\paragraph{Stability of the analysis covariance.} Let $\|\cdot\|$ denote the induced spectral norm for matrices and $\lambda_{\min}(\cdot)$ denote the smallest eigenvalue of a symmetric matrix. For $n>\ell\in\mathbb{N}_0$, set $A_{n,\ell}\coloneqq A_nA_{n-1}\cdots A_{\ell+1}$, $A_{\ell,n}\coloneqq A_{n,\ell}^{-1}$ or invertible $A_n$, and $A_{n,n}\coloneqq \mathbf{I}_{d_X\times d_X}$. For a window length $m\in\mathbb{N}_0$ and $n\geq m+1$,
\begin{equation*}
  O_{n,n-m}\coloneqq\sum_{\ell=n-m}^{n}(H_\ell A_{\ell,n})^TR_\ell^{-1}H_\ell A_{\ell,n}
\end{equation*}
denotes the observability Gramian.

\begin{assumption}{5.O}[Uniform observability]\label{ass: control}
The system in \cref{ass: LGSS} is such that
\begin{itemize}
  \item $A_n$ is invertible and there is $s_{\max}<\infty$ with $\|A_n^{-1}\|^2\|Q_n\|\leq s_{\max}$ for all $n\in\mathbb{N}$,
  \item there are $m\in\mathbb{N}_0$ and $\lambda_O>0$ with $\lambda_{\min}(O_{n,n-m})\geq\lambda_O$ for all $n\geq m+1$, and
  \item there are $\lambda_Q,\lambda_R>0$ and $h_{\max}<\infty$ with $\lambda_{\min}(Q_n)\geq\lambda_Q$, $\lambda_{\min}(R_n)\geq\lambda_R$, and $\|H_n\|\leq h_{\max}$ for all $n\in\mathbb{N}$.
\end{itemize}
\end{assumption}
The first point allows the precision update in \cref{thm: Conjugacy} to be written backwards in time. The second point is the classical uniform observability condition (see e.g. \cite{anderson2012optimal}). The third point, replaces the usual uniform controllability condition by uniformly non-degenerate noise in \cref{ass: LGSS}.

\begin{theorem}[Stability of the DSM covariance]\label{thm: stability}
Given \cref{ass: LGSS} and \cref{ass: control}. Let $k$ be the IMQ kernel in \cref{eqn: IMQ kernel} with $q^2>0$. For the DSM KF in \cref{thm: conjF} with observations $(Y_n)_{n\in\mathbb{N}}$ generated by an arbitrary true DGP $\pi$, let $P^a_n$ denote the analysis covariance and $J^a_n=(P^a_n)^{-1}$ the analysis precision. Then
\begin{enumerate}[label=(\roman*)]
  \item $\|J^a_n\|$ is uniformly bounded for all $n\in\mathbb{N}$ by a deterministic constant.
  \item If $\sup_{n\in\mathbb{N}}\mathbb{E}\left[\|Y_n-H_nm^f_n\|_2^2\right]<\infty$, then $\sup_{n\geq m+1}\mathbb{E}\left[\|P^a_n\|\right]<\infty$ and $(P^a_n)_{n\in\mathbb{N}}$ is weakly stochastically bounded.
\end{enumerate}
\end{theorem}

\begin{proof}
For every realisation, $J^a_n$ follows the regular KF recursion with $R_n$ replaced by $N_n(Y_n)$. Let $\Gamma_n\coloneqq H_n^TN_n^{-1}(Y_n)H_n$ and $S_n\coloneqq A_n^{-1}Q_n A_n^{-T}$. Recall, via the information form analysis step in \cref{thm: Conjugacy},
\begin{equation*}
    J^a_n=(P^f_n)^{-1}+H_n^TN_n^{-1}(Y_n)H_n= A_n^{-T}\left(P^a_{n-1}+S_n\right)^{-1}A_n^{-1}+\Gamma_n\,.
\end{equation*}
Since $\lambda_{\min}(A_n P^a_{n-1}A_n^T+Q_n)\geq\lambda_{\min}(Q_n)\geq\lambda_Q$ and $\|\Gamma_n\|\leq 2k_n^2(Y_n)\,h_{\max}^2\|R_n^{-1}\|\leq 2h_{\max}^2\lambda_R^{-1}$, it follows that $\|J^a_n\|\leq j_{\max}\coloneqq\lambda_Q^{-1}+2h_{\max}^2\lambda_R^{-1}$ for all $n\in\mathbb{N}$ and thus (i). 

For the second expression and $n\geq2$, observe 
\begin{equation*}
    P^a_{n-1}+S_n = (P^a_{n-1})^{1/2}(\mathbf{I}_{d_X\times d_X}+T_n)(P^a_{n-1})^{1/2}
\end{equation*}
with $T_n\coloneqq(J^a_{n-1})^{1/2}S_n(J^a_{n-1})^{1/2}$ positive semi-definite, $\|T_n\|\leq s_{\max}j_{\max}$, and 
\begin{equation*}
    \lambda_{\min}\left((\mathbf{I}_{d_X\times d_X}+T_n)^{-1}\right)\geq\lambda_T\coloneqq(1+s_{\max}j_{\max})^{-1}\,.    
\end{equation*}
Hence, $v^T(P^a_{n-1}+S_n)^{-1}v\geq\lambda_T\,v^TJ^a_{n-1}v$ for all $v\in\mathbb{R}^{d_X}$. Let 
\begin{equation*}
    \widetilde O_{n,n-m}\coloneqq\sum_{\ell=n-m}^{n}A_{\ell,n}^T\Gamma_\ell A_{\ell,n} = \sum_{\ell=n-m}^{n}2k_\ell^2(Y_\ell)\,(H_\ell A_{\ell,n})^TR_\ell^{-1}H_\ell A_{\ell,n}
\end{equation*}
denote the DSM observability Gramian. For $n\geq m+1$ and $u\in\mathbb{R}^{d_X}$,
\begin{align*}
    u^TJ^a_nu&=(A_n^{-1}u)^T(P^a_{n-1}+S_n)^{-1}A_n^{-1}u+u^T\Gamma_nu\geq\lambda_T\,(A_n^{-1}u)^TJ^a_{n-1}A_n^{-1}u+u^T\Gamma_nu\\
    &\geq\lambda_T^2\,(A_{n-1}^{-1}A_n^{-1}u)^TJ^a_{n-2}A_{n-1}^{-1}A_n^{-1}u+\lambda_T(A_{n}^{-1}u)^T\Gamma_{n-1}A_{n}^{-1}u+u^T\Gamma_nu\\
    &\geq\ldots\geq\lambda_T^m\,u^T\widetilde O_{n,n-m}u
\end{align*}
via $\lambda_T\leq1$ and ending in $u^TJ^a_{n-m}u\geq u^T\Gamma_{n-m}u$.

Observe, $\lambda_{\min}(\widetilde O_{n,n-m})\geq 2\lambda_O\min_{n-m\leq\ell\leq n}k_\ell^2(Y_\ell)$ and via $\lambda_{\min}(\Sigma_\ell)\geq\lambda_{\min}(R_\ell)\geq\lambda_R$, $k_\ell^{-2}(Y_\ell)\leq 1+\frac{1}{q^2\lambda_R}\|Y_\ell-H_\ell m^f_\ell\|_2^2$ follows
\begin{equation*}
    \|P^a_n\|=\lambda_{\min}(J^a_n)^{-1}\leq \frac{1}{2\lambda_T^m\lambda_O}\left(1+\frac{1}{q^2\lambda_R}\sum_{\ell=n-m}^{n}\|Y_\ell-H_\ell m^f_\ell\|_2^2\right)\,.
\end{equation*}

Furthermore,
\begin{equation*}
    \sup_{n\geq m+1}\mathbb{E}\left[\|P^a_n\|\right]\leq\frac{1}{2\lambda_T^m\lambda_O}\left(1+\frac{m+1}{q^2\lambda_R}\sup_{n\in\mathbb{N}}\mathbb{E}\left[\|Y_n-H_nm^f_n\|_2^2\right]\right)<\infty\,.
\end{equation*}
Weak stochastic boundedness of $(P^a_n)_{n\geq m+1}$ follows via Markov's inequality. The finitely many $P^a_\ell$ with $\ell\leq m$ do not affect weak stochastic boundedness.
\end{proof}

A union bound over the observability window shows that weakly stochastically bounded innovations $(Y_n-H_nm^f_n)_{n\in\mathbb{N}}$ suffice for $(P^a_n)_{n\in\mathbb{N}}$ to be weakly stochastically bounded. Finite second moments of the innovations are only needed for the bound in expectation. For the regular KF, the chain of arguments results in the deterministic bound $\|P^a_n\|\leq(\lambda_T^m\lambda_O)^{-1}$. For $q^2\to\infty$, so dropping the robustness, the DSM KF recovers the spectral bound of a regular KF with $R_n$ rescaled by $\frac{1}{2}$. Boundedness of the innovations depends on the DSM analysis mean investigated next.

\paragraph{Bounded DSM mean increment.} Given the IMQ kernel in \cref{eqn: IMQ kernel}, observe for the shifted observation,
\begin{equation}\label{eqn: corr obs}
    \tilde{y}-Hm^f=y-2N(y)\nabla_{y}k^2(y)-Hm^f= \left(\mathbf{I}_{d_Y\times d_Y}+\tfrac{2}{q^2}\,k^2(y)\,R\,\Sigma^{-1}\right)(y-Hm^f)\,.
\end{equation}

\begin{lemma}[Bounded influence]\label{lem: influence}
Given \cref{ass: LIP} and let $k$ be the IMQ kernel in \cref{eqn: IMQ kernel} with $q^2>0$. The DSM analysis mean increment $m^a-m^f$ is bounded uniformly in $y\in\mathbb{R}^{d_Y}$ and $\|m^a-m^f\|_2\to 0$ as $\|y-Hm^f\|_2\to\infty$.
\end{lemma}

\begin{proof}
Via $\widetilde{K}(y)=2k^2(y)P^aH^TR^{-1}$ and \cref{eqn: corr obs},
\begin{equation*}
    m^a-m^f=-\widetilde{K}(y)\left[Hm^f-\tilde{y}\right]= 2k^2(y)\,P^aH^TR^{-1}\left(\mathbf{I}_{d_Y\times d_Y}+\tfrac{2}{q^2}k^2(y)R\Sigma^{-1}\right)(y-Hm^f)\,.
\end{equation*}
Since $\|P^a\|\leq\|P^f\|$ and $\|\Sigma^{-1}\|\leq\|R^{-1}\|$ in the single analysis step, and $k^2\leq 1$, all matrix factors are bounded uniformly in $y$. Further,
\begin{equation*}
    k^2(y)\|y-Hm^f\|_2\leq\frac{q^2\|\Sigma\|\|y-Hm^f\|_2}{q^2\|\Sigma\|+\|y-Hm^f\|_2^2}\,,
\end{equation*}
which is bounded in $y$ and tends to $0$ as $\|y-Hm^f\|_2\to\infty$.

\end{proof}

The regular KF mean increment $K(y-Hm^f)$ is unbounded in $y$. \Cref{lem: influence} is the filtering counterpart to \cref{thm: bias rob} for the analysis mean. Given \cref{ass: control} and additionally $\sup_{n\in\mathbb{N}}\|A_n\|<1$ and $\sup_{n\in\mathbb{N}}\|R_n\|<\infty$, the mean increments are bounded uniformly in $n$ and $y$, and the analysis mean remains bounded for arbitrary true DGP $\pi$. The innovation condition in \cref{thm: stability} then holds whenever $\sup_{n\in\mathbb{N}}\mathbb{E}\left[\|Y_n\|_2^2\right]<\infty$. Generalising the results to an analysis of the DSM error dynamics needs controlling the innovation. This requires specifying a model for the true DGP $\pi$ and is beyond the aims here (see \cref{sec: discussion}).

\subsection{Tuning and behaviour for high-dimensional observations}\label{sec: tuning}

The DSM KF assumes only limited knowledge of the observation error. This opens an additional degree of freedom for adjusting information processing, here reduced to the scalar $q^2$ in \cref{eqn: IMQ kernel}. \Cref{thm: bias rob} holds for any finite $0<q^2$, yet the choice affects accuracy. We investigate influence as well as propose a heuristic default choice of the parameter $q^2$.

The parameter appears in similar form in \cite{altamirano2023robust,altamirano2023robust2,duran2024outlier} with tuning only discussed to a minor extent. We follow the line of thought in \cite{altamirano2023robust} by choosing $q^2$ with respect to the well-specified case. The Mahalanobis distance in the IMQ weight function serves as a whitening of the observation under the well-specified observation marginal $Y_n\sim \mathcal{N}(H_nm^f_n,\Sigma_n)$ with respect to the true DGP $\pi$ for a single analysis step. Then $\Theta_n\coloneqq\|Y_n-H_nm^f_n\|^2_{\Sigma_n^{-1}}\sim\chi^2(d_Y)$ with $\mathbb{E}_{Y_n}\left[\Theta_n\right]=d_Y$ and $\mathrm{Var}_{Y_n}\left(\Theta_n\right)=2d_Y$.

We investigate the effect of the parameter choice on the precision update in \cref{thm: Conjugacy} in expectation with respect to the observation marginal, so
\begin{equation}\label{eqn: eprec}
    \mathbb{E}_{Y_n}\left[J^a_n\right]=\mathbb{E}_{Y_n}\left[J^f_n+H^T_nN_n^{-1}(Y_n)H_n\right]=J^f_n+H^T_n\mathbb{E}_{Y_n}\left[2k_n^2(Y_n)\right]R_n^{-1}H_n.
\end{equation}
We simplify notation replacing $2k^2_n(Y_n)$ by $g(\Theta)=2\left(1+\frac{\Theta}{q^2}\right)^{-1}$ with $\Theta\sim\chi^2(d_Y)$ and support $\vartheta\in[0,\infty)$. The function $g$ is convex, $g(\vartheta)\in(0,2]$, and $g(\mu)=2\left(1+\frac{d_Y}{q^2}\right)^{-1}$ for $\mu\coloneqq\mathbb{E}\left[\Theta\right]$. Via Jensen's inequality, 
\begin{equation*}
    g(\mu)=2\left(1+\frac{d_Y}{q^2}\right)^{-1}\leq\mathbb{E}\left[g(\Theta)\right]
\end{equation*}
showing that $\mathbb{E}[g(\Theta)]\geq1$ for $q^2\geq d_Y$, i.e. at least the precision increase of the regular KF in expectation. This indicates overconfidence compared to the (information optimal) precision update of the regular KF. To obtain an upper bound, we consider the Jensen gap $G=\mathbb{E}\left[g(\Theta)\right]-g(\mu)$ and introduce a technical lemma based on \cite{gao2017bounds}.

\begin{lemma}[Upper bound on the Jensen gap]\label{thm: JGap1}
    Let $g(z)=2\left(1+\frac{z}{q^2}\right)^{-1}$ with $q^2>0$. Let $Z$ be some non-negative random variable with mean $\mathbb{E}[Z]=\mu<\infty$ and variance $\mathrm{Var}(Z)=\sigma^2<\infty$. The Jensen gap $G=\mathbb{E}[g(Z)]-g(\mu)$ can be upper bounded via
    \begin{equation*}
        G\leq\frac{2\sigma^2}{(q^2+\mu)^2}
    \end{equation*}
    leading to the two-sided Jensen-type inequality
    \begin{equation*}
        g(\mu)\leq\mathbb{E}[g(Z)]\leq g(\mu)+\frac{2\sigma^2}{(q^2+\mu)^2}.
    \end{equation*}
\end{lemma}

\begin{proof}
As $\mathbb{E}[Z-\mu]=0$ and $g'(\mu)=-\frac{2}{q^2}(1+\frac{\mu}{q^2})^{-2}$,
\begin{align*}
    G&=\mathbb{E}\left[g(Z)-g(\mu)-g'(\mu)(Z-\mu)\right]=2q^2\mathbb{E}\left[\frac{(Z-\mu)}{(q^2+\mu)^2}-\frac{(Z-\mu)}{(q^2+Z)(q^2+\mu)}\right]\\
    &=2q^2\mathbb{E}\left[\frac{(Z-\mu)^2}{(q^2+\mu)^2(q^2+Z)}\right]\leq\frac{2\sigma^2}{(q^2+\mu)^2}
\end{align*}
via $q^2+Z\geq q^2$ (see \cite[eq. (6)]{gao2017bounds}).
\end{proof}

The DSM KF precision update in the well-specified setting then satisfies
\begin{equation*}
    2\left(1+\frac{d_Y}{q^2}\right)^{-1}\leq\mathbb{E}_{Y_n}\left[2k^2_n(Y_n)\right]\leq2\left(1+\frac{d_Y}{q^2}\right)^{-1}+\frac{4d_Y}{(q^2+d_Y)^2}\,.
\end{equation*}
This suggests choosing the parameter $q^2$ proportional to the observation dimension $d_Y$ with the choice $q^2=d_Y$ considered in the following statement.

\begin{proposition}[Analysis precision for well-specified observation likelihoods]\label{thm: hiDim}
Assume for a single analysis step of \cref{thm: conjF}, that the observation marginal, i.e. $Y_n\sim \mathcal{N}(H_nm^f_n,\Sigma_n)$, is well-specified. Further, let the weight kernel $k$ be as in \cref{eqn: IMQ kernel} with $q^2=d_Y$. Then
\begin{equation*}
    \mathbb{E}_{Y_n}\left[J^a_n\right]=J^f_n+\rho_{d_Y}H_n^TR_n^{-1}H_n\quad\text{with}\quad \rho_{d_Y}=\mathbb{E}_{Y_n}\left[2k_n^2(Y_n)\right]\quad\text{and}\quad 1\leq\rho_{d_Y}\leq1+\frac{1}{d_Y}\,.
\end{equation*}
Therefore, the single analysis step precision update of the DSM KF is asymptotically unbiased with respect to the regular Kalman filter as $d_Y\to\infty$ granted the assumptions.
\end{proposition}

The proof follows from the previous considerations in \cref{eqn: eprec} and \cref{thm: JGap1}. Overconfidence in the precision update for the default choice $q^2=d_Y$ in the well-specified case is especially present for $d_Y$ small with $\rho_1\approx1.31$. Lowering $q^2$ to reduce overconfidence in precision does not increase accuracy due to the factor $\frac{2}{q^2}$ in the shifted observation in \cref{eqn: corr obs}. The DSM KF should mainly be considered over the regular KF when there is concern for misspecification. We then expect the down-weighting of outliers and the overconfidence of the default choice $q^2=d_Y$ to balance.

\section{Ensemble approximations of the analysis step}\label{sec: ens nl}

The family of ensemble Kalman filters (EnKFs) was popularised by its empirical success via Gaussian approximation. Departing from \cref{ass: LGSS} and assuming non-linear signal dynamics, the forecast distribution is generally no longer Gaussian. Accordingly, the forecast covariance matrix is no longer available and propagating from analysis to forecast via linearisation becomes costly in high dimensions. EnKFs instead propagate the previous analysis ensemble according to the signal dynamics to obtain a forecast ensemble and access empirical moments. The DSM analysis step in \cref{thm: conjF}, just as its regular counterpart, only requires access to the forecast mean and covariance matrix. The resulting analysis step is subject to both empirical and Gaussian approximation and only remains consistent in the ensemble limit under linear signal dynamics.

We consider the stochastic EnKF (also EnKF with perturbed observations) and one type of deterministic EnKF, the ensemble square root filter (ESRF) (see e.g. \cite{reich2015probabilistic}), and its modification for analysis in ensemble subspace, basis of the LETKF (see \cite{hunt2007efficient}). We only consider linear observation operators in the following. Non-linear observation operators may be considered via the usual modifications employing cross-covariances (see e.g. \cite{evensen2022data}) or local linear approximation (see \cite{hunt2007efficient}).

\subsection{Stochastic and deterministic ensemble coupling}

We follow the derivation in \cite{reich2015probabilistic} via coupling forecast and analysis ensemble either in a stochastic or deterministic sense to adjust the DSM Kalman filter parallel to its regular counterparts. For stochastic coupling, assume access to a forecast ensemble $\{x^{f,(i)}\}_{i=1}^M$ with ensemble size $M$ and notation dropping time indices and the empirical forecast mean $\bar{x}_M^f$ and empirical covariance matrix $P_M^f$. We update the forecast ensemble via
\begin{equation*}
   x^{a,(i)}=x^{f,(i)}-\widetilde{K}(y)[Hx^{f,(i)}+\xi^{(i)}-\tilde{y}]
\end{equation*}
with $\{\xi^{(i)}\}_{i=1}^M$ \emph{i.i.d.} draws from $\Xi\sim\mathcal{N}(0,N(y))$ and $P^f$ replaced by $P_M^f$. The ensemble update is consistent with respect to the DSM KF given \cref{ass: LGSS},
\begin{align*}
    &\mathbb{E}\left[X^a\right]=m^f-\widetilde{K}(y)\left[Hm^f-\tilde{y}\right]\ \text{and}\\
    &\mathbb{E}\left[\left(X^a-\mathbb{E}[X^a]\right)\left(X^a-\mathbb{E}[X^a]\right)^T\right]=P^f-\widetilde{K}(y)HP^f\,.
\end{align*}
Mean consistency directly follows from linearity and covariance consistency exploits independence of $\Xi$ and $X^f$ plus a matrix identity. We refer to \cite{reich2015probabilistic} for additional details. Note that this strictly requires the weight kernel $k$ to be non-stochastic with respect to the forecast ensemble. In practice, the IMQ weight kernel $k$ needs adjusting via replacing the forecast parameters with their empirical counterparts,
\begin{equation*}
    k_M^2(y)=\left(1+\frac{\|y-H\bar{x}_M^f\|_{\Sigma_M^{-1}}^2}{q^2}\right)^{-1}\quad\text{with}\quad\Sigma_M=HP^f_MH^T+R\,.
\end{equation*}
The adjusted Kalman gain $\widetilde{K}(y)$ and the shifted observation $\tilde{y}$ are then quantities strongly depending on the forecast ensemble and consistency holds only in the ensemble limit. Similarly, the tuning heuristic does then not directly translate as the finite sample empirical counterparts induce additional variance. Nevertheless, we propose following the default choice of $q^2=d_Y$.

The stochastic DSM EnKF stated here is referred to as \emph{average-particle} EnKF in \cite{duran2024outlier} indicating the decision to replace the forecast mean in the weight kernel by the empirical forecast mean. In contrast, the \emph{per-particle} EnKF replaces the forecast mean in the weight kernel with the individual ensemble members.The \emph{average-particle} EnKF is theoretically more consistent with the marginal standardisation $\|\cdot\|^2_{\Sigma_n^{-1}}$ employed in the weight kernel. The \emph{per-particle} EnKF requires a conditional standardisation $\|\cdot\|^2_{R_n^{-1}}$. The respective argument arises via the heuristic of measuring deviation beyond the well-specified case either in $p(y_n)$ or $p(y_n|x^f_n)$. The conditional standardisation does still suit the tuning idea in \cref{sec: tuning}. Under non-linear dynamics and observation operators, the \emph{average-particle} standardisation relies on an empirical, Gaussian approximation of the observation marginal. The \emph{per-particle} standardisation only considers departure from Gaussian observation noise.

While both options are available for the stochastic EnKF, we keep the \emph{average-particle} approach for constructing the deterministic DSM ensemble square root filter. The stochastic EnKF is known to be less accurate in high dimension and for small ensemble sizes. The ESRF maintains accuracy even for comparatively small ensemble sizes (see e.g. \cite{reich2015probabilistic}).

We consider the DSM analysis parameters given the empirical forecast mean and covariance matrix with
\begin{itemize}
    \item adjusted Kalman gain $\widetilde{K}(y)=P^f_MH^T\left[N(y)+HP^f_MH^T\right]^{-1}$,
    \item analysis covariance $\widetilde{P}^a=P^f_M-\widetilde{K}(y)HP^f_M$, and
    \item analysis mean $\widetilde{m}^a=\bar{x}_M^f-\widetilde{K}(y)[H\bar{x}_M^f-\tilde{y}]$.
\end{itemize}
 Let $X^f\in\mathbb{R}^{d_X\times M}$ be the matrix of forecast ensemble anomalies, so $X^f$ has columns $(X^f)_i=x^{f,(i)}-\bar{x}_M^f$. We are interested in the linear transformation of the ensemble via
\begin{equation*}
    x^{a,(i)}=\widetilde{m}^{a}+(X^fS)_i\quad\text{with}\quad S\ \text{s.t.}\ \widetilde{P}^a=\frac{1}{M-1}X^fS(X^fS)^T\,.
\end{equation*}
Exploiting $P^f_M=\frac{1}{M-1}X^f(X^f)^T$, this holds for $SS^T=T$ with
\begin{equation*}
    T\coloneqq \mathbf{I}_{M\times M}-\frac{1}{M-1}(HX^f)^T\left[HP^f_MH^T+N(y)\right]^{-1}HX^f\,.
\end{equation*}
We take the unique symmetric positive definite square root $S=T^{\frac{1}{2}}$ keeping the analysis mean $\widetilde{m}^a$. Note, adapting the ESRF is only implicit in the adjusted parameter updates. Again, we refer to \cite{reich2015probabilistic} for additional details.

\subsection{Analysis in ensemble subspace}\label{sec: DSM LETKF}

As a common extension to the ESRF, we consider analysis in ensemble anomaly subspace. Localisation is not discussed here, yet practice directly transfers and we refer to \cite{hunt2007efficient} for details. The approach starts with an optimisation-centric view of the DSM posterior in the Bayesian inverse problem. Consider the DSM negative log-posterior (proof of \cref{thm: Conjugacy})
\begin{equation*}
    \mathcal{L}_\mathrm{DSM}(x)\coloneqq\frac{1}{2}(x-\bar{x}_M^f)^T(P_M^f)^{-1}(x-\bar{x}_M^f)+\frac{1}{2}(y-Hx)^TN(y)^{-1}(y-Hx)+(Hx)^T\nabla_y2k_M^2(y)\,.
\end{equation*}
Via orthogonal projection into ensemble subspace the problem is transferred to the objective
\begin{equation*}
    \mathcal{L}^*_\mathrm{DSM}(v)\coloneqq\frac{1}{2}(M-1)v^Tv+\frac{1}{2}([y-\bar{y}^f_M]-Y^fv)^TN(y)^{-1}([y-\bar{y}^f_M]-Y^fv)+(Y^fv)^T\nabla_y2k_M^2(y)
\end{equation*}
with $v\in\mathbb{R}^M$, $\bar{y}_M^f=H\bar{x}_M^f$ and $Y^f=HX^f$. The objective $\mathcal{L}^*_\mathrm{DSM}$ is solved by the DSM analysis step via
\begin{itemize}
    \item $\widetilde{K}_v(y)=\frac{1}{M-1}(Y^f)^T\left[N(y)^{-1}+\frac{1}{M-1}Y^f(Y^f)^T\right]^{-1}$,
    \item $\widetilde{P}_v^a=\frac{1}{M-1}\mathbf{I}_{M\times M}-\frac{1}{M-1}\widetilde{K}_v(y)Y^f$, and
    \item $\widetilde{m}^a_v=\widetilde{K}_v(y)\left[(y-\bar{y}^f_M)-2N(y)\nabla_{y}k_M^2(y)\right]$
\end{itemize}
with counterparts in signal space given via $\widetilde{P}^a=X^f\widetilde{P}_v^a(X^f)^T$ and $\widetilde{m}^a=\bar{x}_M^f+X^f\widetilde{m}^a_v$. The analysis ensemble is then obtained e.g. via the ESRF utilising the obtained analysis parameters.

In addition to the DSM variants, the WoLF variants of the KF and EnKF in \cite{duran2024outlier} are investigated in the next section. A novel WoLF variant of the LETKF constructed via similar arguments is considered. 

\section{Simulation experiments}\label{sec: simulation study}

The simulation experiments are not an exhaustive comparative study, but serve to convey dynamics and challenges of the presented approaches with a focus on GBI-based EnKFs for non-linear dynamical systems.

\paragraph{Experimental design.}
We consider the regular, WoLF and DSM variants of the KF, EnKF, and LETKF. The simulation study includes linear dynamics via 2D target tracking and non-linear signal processes via the stochastic Lorenz-63 system to assess the performance of ensemble approximations exposed to chaotic signal dynamics and the Lorenz-96 system with $d_X=d_Y=40$ to assess the impact of higher dimension. Each experiment covers different severities of contamination. For the KF and LETKF variants, we additionally investigate the empirical coverage of credible intervals, and for the EnKF the effect of ensemble size.

Unless stated otherwise, the DSM variants employ the default tuning $q^2=d_Y$ (\cref{sec: tuning}), and the WoLF variants, implemented as in \cite{duran2024outlier}, employ the threshold $c^2=d_Y$.

To capture uncertainty quantification (UQ) via the Gaussian approximation, we mainly consider the component-average continuous ranked probability score (CRPS; see e.g. \cite{gneiting2007strictly}),
\begin{equation*}
    \mathrm{CRPS}(x^\dagger;\, m^a,P^a)=\frac{1}{d_X}\sum_{i=1}^{d_X}\sigma_i\left[ z^\dagger_i\left(2\Phi\left( z^\dagger_i\right)-1\right)+2\phi\left(z^\dagger_i\right)-\frac{1}{\sqrt{\pi}}\right]
\end{equation*}
with $z^\dagger_i\coloneqq\frac{x^\dagger_i-m^a_i}{\sqrt{\sigma^2_i}}$ for true state $x^\dagger$, analysis mean $m^a$, analysis variances $\sigma_i^2=P^a_{ii}$, and $\phi$ and $\Phi$ the standard Gaussian density and CDF. For the EnKF and LETKF, mean and covariance are replaced by their empirical counterparts. The CRPS is averaged over the analysis times and summarised by the median over repetitions. Note that the median is taken over all runs, including those that lose track of the signal. As an alternative to the Gaussian CRPS, the unbiased ensemble CRPS in \cite{ferro2014fair} can be considered for the EnKF and LETKF, leading to the same conclusions. The time-averaged root mean squared error (RMSE) is additionally employed for evaluation of the ensemble size comparison of the EnKF. Additional evaluation metrics are provided in the Supplementary Material \cite{supplement}.

Misspecification of the observation noise is simulated via contamination of the Gaussian observation error by an additional, highly inflated Gaussian error. Severity of misspecification is steered by frequency $\varepsilon\in[0,1]$ and degree of inflation $\lambda\geq1$,
\begin{equation*}
    Y^{\varepsilon,\lambda}_n=H_nX_n+R^{\frac{1}{2}}_nV_n^{\varepsilon,\lambda}\quad\text{with}\quad
    V_n^{\varepsilon,\lambda}\sim (1-\varepsilon)\ \mathcal{N}(0,\mathbf{I}_{d_Y\times d_Y})+\varepsilon\ \mathcal{N}(0,\lambda\cdot\mathbf{I}_{d_Y\times d_Y}).
\end{equation*}

\subsection{Observation noise misspecification in linear filtering}\label{sec: simulation study LGSS}

\begin{figure}[h!]
    \centering
    \includegraphics[width=\linewidth]{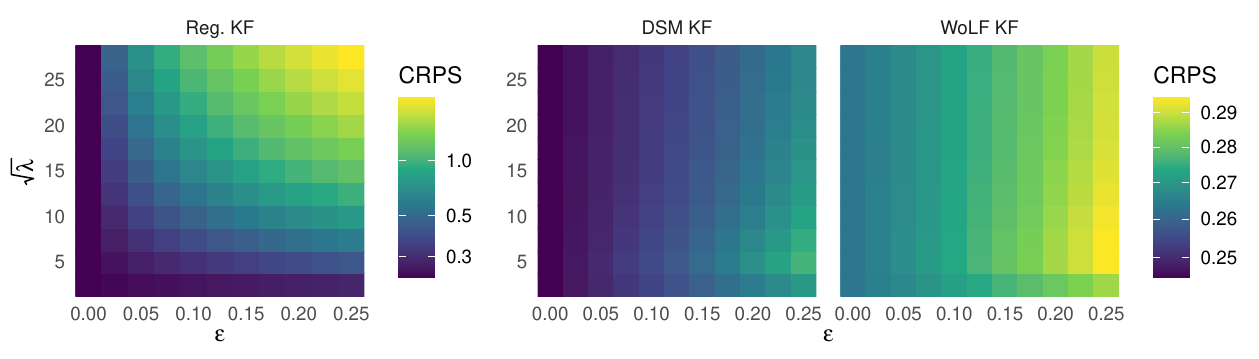}
    \caption{Median CRPS over $M_\mathrm{MC}=2500$ repetitions of the regular KF (left), DSM KF (middle), and WoLF KF (right) for different frequencies ($\varepsilon$) and degrees ($\lambda$). The DSM and WoLF KF share the right scale.}
    \label{fig: TT IC}
\end{figure}

The main interest given \cref{ass: LGSS} lies in capabilities for UQ under contamination and default tuning choice. The 2D target tracking example is adapted from \cite{sarkka2023bayesian}. The signal is discretised with $T_{\mathrm{end}} = 50$ and discretisation time step $\Delta t = 0.1$. The model components are
\begin{equation*}
    A=\begin{bmatrix}1 & 0 & \Delta t & 0 \\0 & 1 & 0 & \Delta t \\0 & 0 & 1 & 0 \\0 & 0 & 0 & 1\end{bmatrix},\quad Q=\begin{bmatrix}\frac{\Delta t^3}{3} & 0 & \frac{\Delta t^2}{2} & 0 \\
0 & \frac{\Delta t^3}{3} & 0 & \frac{\Delta t^2}{2} \\\frac{\Delta t^2}{2} & 0 & \Delta t & 0 \\
0 & \frac{\Delta t^2}{2} & 0 & \Delta t\end{bmatrix},\quad\text{and}\quad x_0=\begin{bmatrix}0 \\0 \\1 \\1\end{bmatrix}
\end{equation*}
and observation components 
\begin{equation*}
    H=\begin{bmatrix}1 & 0 & 0 & 0 \\0 & 1 & 0 & 0
\end{bmatrix},\quad\text{and}\quad R= \begin{bmatrix}0.1 & 0.1^2 \\0.1^2 & 0.1\end{bmatrix}\,.
\end{equation*}
Observations are produced after each time step. The KF variants are initialised with parameters $m_0=x_0$ and $P_0= \mathbf{I}_{4\times4}$. The experiment was repeated $M_\mathrm{MC}=2500$ times for different severities $\varepsilon\in\{0, 0.025, 0.05, \ldots, 0.25\}$ and $\sqrt{\lambda}\in\{2.5, 5, 7.5, \ldots, 27.5\}$. The results are presented for the CRPS in \cref{fig: TT IC} and for coverage in \cref{tab: TT coverage}.

\begin{table}[ht]
\centering
\begin{tabular}{lcc}
\toprule
 & well-specified & all \\
\midrule
KF & 0.90 & 0.64 (0.27--0.90) \\
DSM KF & 0.86 & 0.86 (0.84--0.86) \\
WoLF KF & 0.92 & 0.92 (0.90--0.92) \\
\bottomrule
\end{tabular}
\caption{Empirical coverage of the componentwise $90\%$ credible intervals in the well-specified case and over all severities in \cref{fig: TT IC} (median and range).}
\label{tab: TT coverage}
\end{table}

We observe the robustness of the GBI-based methods with the DSM KF slightly but consistently better in Gaussian CRPS. Regarding UQ as indicated by \cref{tab: TT coverage}, the DSM KF is overconfident to a minor degree whereas the WoLF KF coverage is more conservative. Both methods perform worst for highly frequent but moderate degree contamination. Scaling the WoLF observation log-likelihood by a learning rate of $2$ (see e.g. \cite{altamirano2023robust}) likely recovers the minor difference observed in CRPS, yet likely also leads to overconfidence. A relevant impact of the shifted observation in the DSM KF cannot be reliably deduced from the results.

\subsection{Observation noise misspecification in non-linear filtering}\label{sec: simulation study NL}

\subsubsection{Stochastic Lorenz-63 model}

We adjust the model setup in \cite{reich2015probabilistic} with Euler-Maruyama discretization. Trajectories are taken over time windows of length $T_{\mathrm{end}} = 50$ with discretisation step size $\Delta t = 0.001$ and observations generated every $t_{\mathrm{out}}=0.05$, so every $50$ steps. The signal is updated according to 
\begin{equation*}
    x_{n} = x_{n-1} + \Delta t f(x_{n-1}) + \sqrt{\Delta t} W_n\quad\text{with}\quad f(x) = \begin{bmatrix} 10(x_2 - x_1) \\x_1(28 - x_3) - x_2 \\x_1x_2 - (8/3)x_3\end{bmatrix},
\end{equation*}
$W_n\overset{\mathit{i.i.d.}}{\sim}\mathcal{N}(0,\mathbf{I}_{3\times3})$, and $x_0 = (-0.587, -0.563, 16.87)^T$. Observations are produced with 
\begin{equation*}
    H=\begin{bmatrix} 1 & 0 & 0 \end{bmatrix}\quad \text{and}\quad R = 0.5
\end{equation*}
up to contamination. We compare the regular, DSM, and WoLF variants of the stochastic EnKF with an ensemble size of $M_\mathrm{ens}=10$. We employ the \emph{average-particle} variants of the DSM and WoLF EnKF. The initial ensemble is produced by sampling $x_0^{a,(i)}\overset{\mathit{i.i.d.}}{\sim}\mathcal{N}(x_0,0.1\cdot\mathbf{I}_{3\times 3})$ for $i=1,2,\ldots,M_\mathrm{ens}$. No covariance inflation is employed and the first $10\%$ of the analysis cycles are excluded from evaluation as spin-up.

\begin{figure}[h!]
    \centering
    \includegraphics[width=1\linewidth]{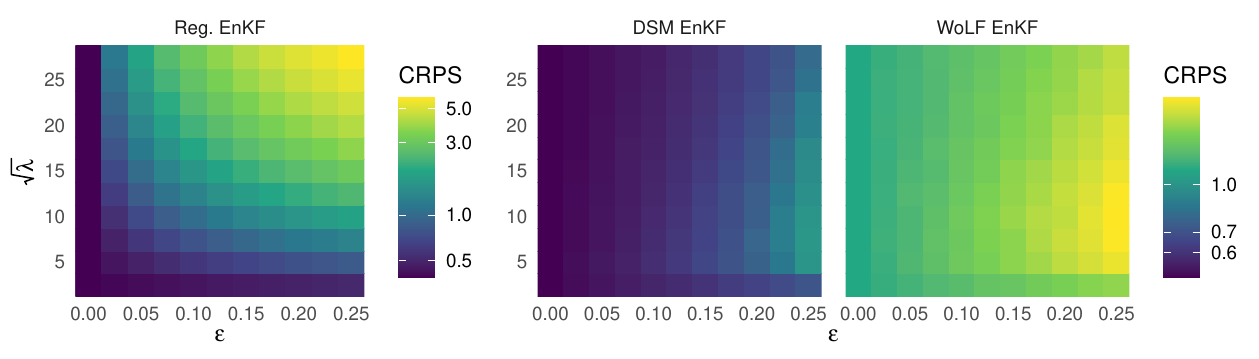}
    \caption{Median CRPS over $M_\mathrm{MC}=1000$ repetitions of the stochastic variants of the regular EnKF (left), DSM EnKF (middle), and WoLF EnKF (right) for different frequencies ($\varepsilon$) and degrees ($\lambda$). The DSM and WoLF EnKF share the right scale.}
    \label{fig: Lor63 IC}
\end{figure}

\begin{figure}[h!]
    \centering
    \includegraphics[width=1\linewidth]{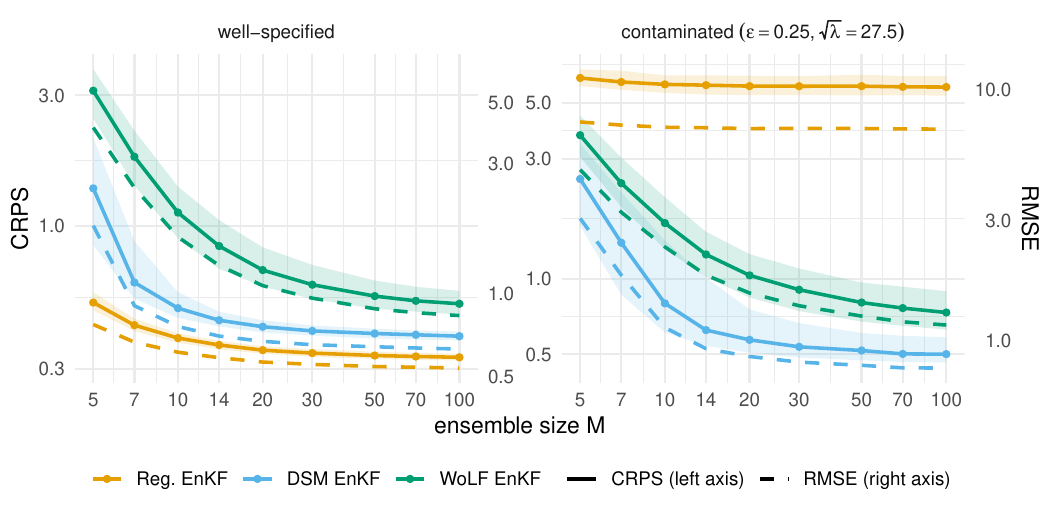}
    \caption{Median CRPS (solid, shaded $25$--$75\%$ quantile band; left $y$-axis) and RMSE (dashed; right $y$-axis) over $M_\mathrm{MC}=1000$ repetitions of stochastic EnKF variants for different ensemble sizes in the well-specified (left) and contaminated model (right).}
    \label{fig: ens EnKF IC}
\end{figure}

The experiment was repeated $M_\mathrm{MC}=1000$ times for different severities $\varepsilon\in\{0, 0.025, 0.05, \ldots, 0.25\}$ and $\sqrt{\lambda}\in\{2.5, 5, 7.5, \ldots, 27.5\}$ including the well-specified case for $\varepsilon=0$. The results are presented for the CRPS in \cref{fig: Lor63 IC}. Similarly, the experiment was repeated $M_\mathrm{MC}=1000$ times for different ensemble sizes $M_\mathrm{ens}\in\{5, 7, 10, 14, 20, 30, 50, 70, 100\}$ in the well-specified model and the contaminated model with frequency $\varepsilon=0.25$ and degree $\lambda = 27.5^2$. The results are presented in \cref{fig: ens EnKF IC}.

The regular stochastic EnKF performs best in the well-specified case with the DSM EnKF second best. Under contaminated observations, we observe the DSM EnKF ahead of the WoLF EnKF across contamination severities and ensemble sizes. We attribute this difference to the property of the DSM EnKF to not just slow down learning but also increase information gain in regimes with no relevant contamination realisations. This may support stability and reduce filter divergence in the notoriously chaotic Lorenz-63 system. Again, adjusting the learning rate of the WoLF KF may account for this. The insight supports the need to not just reduce information intake but make appropriate use of available information where possible in non-linear systems as discussed in \cref{sec: discussion}.

\subsubsection{40D Lorenz-96 model}

We follow the experimental setup in \cite{sun2022control} and refer to their introduction for a more comprehensive overview. Trajectories are taken over time windows of length $T_{\mathrm{end}} = 73$. The model is run via a fourth-order Runge--Kutta scheme with integration time step $\Delta t=0.01$ applied to the $40$ differential equations
\begin{equation*}
    \frac{\mathrm{d}x_i}{\mathrm{d}t} = (x_{i+1}-x_{i-2})x_{i-1}-x_i+F
\end{equation*}
with $x_i$ the $i$-th state component, the convention $x_0=x_{40},\ x_{-1}=x_{39}$ and $x_{41}=x_1$ (resembling states on a latitude circle), and constant forcing $F=8$. Observations are produced every $t_{\mathrm{out}}=0.05$ time units via  
\begin{equation*}
    H=\mathbf{I}_{40\times40}\ \text{and}\  R=\mathbf{I}_{40\times40}. 
\end{equation*}
The initial state $x_0$ is produced by an initial draw from $\mathcal{N}(0,\mathbf{I}_{40\times 40})$ followed by a burn-in with $T_{\mathrm{burn}} = 12.15$.

We implement both multiplicative covariance inflation and localisation with the parameters given in \cite{sun2022control} and use their choice of $M_\mathrm{ens}=10$. Details are provided in the Supplementary Material \cite{supplement}. The DSM and WoLF weights are evaluated for each local analysis from the inflated forecast covariance and the localised $R$. This requires adjusting the tuning heuristic in \cref{sec: tuning} for $q^2$ via replacing $d_Y$ by $d_Y^{\mathrm{eff}}\approx\sum_j\ell_j$, the effective number of observations, with localisation weight $\ell_j\in(0,1]$ of the $j$-th local observation. For localisation as in \cite{sun2022control}, this suggests $q^2=d_Y^{\mathrm{eff}}\approx13.7$. Otherwise, the numerical implementation of the LETKF variants is adapted directly from \cite{hunt2007efficient}. The initial ensemble is produced via $x_0^{a,(i)}\overset{\mathit{i.i.d.}}{\sim}\mathcal{N}(x_0,\mathbf{I}_{40\times 40})$ for $i=1,2,\ldots,M_\mathrm{ens}$. The first $10\%$ of the analysis cycles are again excluded from evaluation as spin-up.

\begin{figure}[h!]
    \centering
    \includegraphics[width=1\linewidth]{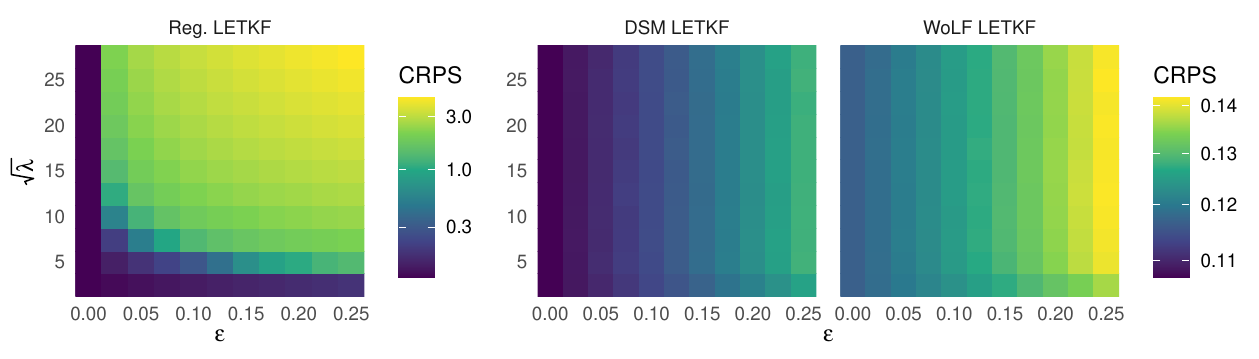}
    \caption{Median CRPS over $M_\mathrm{MC}=500$ repetitions of the regular LETKF (left), DSM LETKF (middle), and WoLF LETKF (right) for different frequencies ($\varepsilon$) and degrees ($\lambda$). The DSM and WoLF LETKF share the right scale.}
    \label{fig: Lor96 IC}
\end{figure}

The experiment was repeated $M_\mathrm{MC}=500$ times for different severities $\varepsilon\in\{0, 0.025,\allowbreak 0.05, \ldots, 0.25\}$ and $\sqrt{\lambda}\in\{2.5, 5, 7.5, \ldots, 27.5\}$. The results are presented for the CRPS in \cref{fig: Lor96 IC} and for coverage in \cref{tab: L96 coverage}.

\begin{table}[ht]
\centering
\begin{tabular}{lcc}
\toprule
 & well-specified & all \\
\midrule
LETKF & 0.92 & 0.12 (0.06--0.92) \\
DSM LETKF & 0.92 & 0.91 (0.90--0.92) \\
WoLF LETKF & 0.98 & 0.98 (0.97--0.98) \\
\bottomrule
\end{tabular}
\caption{Empirical coverage of the componentwise $90\%$ credible intervals in the well-specified case and over all severities in \cref{fig: Lor96 IC} (median and range).}
\label{tab: L96 coverage}
\end{table}

The regular LETKF performs best in the well-specified case with the DSM LETKF a very close second and matching uncertainties. Under contamination, the regular LETKF breaks down and loses track of the signal, with the coverage of its $90\%$ credible intervals collapsing. Both robust variants remain stable across all severities with the DSM LETKF ahead of the WoLF LETKF in Gaussian CRPS. Furthermore, the DSM LETKF maintains reliable coverage whereas the WoLF LETKF is again more conservative. Compared to the experiment in \cref{sec: simulation study LGSS}, we no longer observe overconfidence for the DSM LETKF, likely due to the combination of larger $d_Y^{\mathrm{eff}}$ (see \cref{thm: hiDim}; $\rho_{13.7}\approx1.07$) and the multiplicative covariance inflation from \cite{sun2022control} causing even the regular LETKF to be slightly conservative ($92\%$ coverage in the well-specified case).

\section{Discussion}\label{sec: discussion}

\paragraph{Misspecification and error assumptions.}

Global bias robustness of the DSM KF analysis step in \cref{thm: bias rob} specified misspecification to be understood in the sense of Huber's $\varepsilon$-contamination. Following \cite{mu2023huber}, it treats the assumed DGP, characterised by the modelled observation likelihood, as an accurate description of the true DGP up to influential instances, the contaminations. \Cref{thm: bias rob} guarantees that no single instance of contamination can impact the single analysis step beyond a finite degree. \Cref{lem: influence} carries this over to the analysis mean of the DSM KF, again for arbitrary $\pi$.

\Cref{thm: stability} imposes no restrictions on the true DGP $\pi$ for the analysis precision bound in (i), yet requires assumptions on the innovations $Y_n-H_nm^f_n$ for the analysis covariance bound in (ii). Here, weakly stochastically bounded innovations suffice for a weakly stochastically bounded analysis covariance and uniformly bounded second moments of the innovations guarantee a bound in expectation. Stronger results on stability need a model for the true signal, in turn restricting general applicability. Since $N_n(Y_n)$ depends on the filter state, results on the existence and uniqueness of an invariant measure, such as in \cite{solo1996stability}, do not apply directly. Growing errors increase model-observation mismatch and inflation in the analysis step, yet this also increases $P^f_n$ and hence $\Sigma_n$ in the standardisation in $k$, leading to a trade-off that may recover a limiting distribution for suitable true DGPs $\pi$. Additionally, the shift in \cref{eqn: corr obs} perturbs the analysis step by a term of order $\frac{1}{q^2}$. This suggests an additional stability condition of the form $q^2\geq q_0^2$.

\paragraph{Misspecification and signal non-linearity.}

Dropping \cref{ass: LGSS}, intuitions still transfer to the case of non-linear signal dynamics. In the simulation experiments, the DSM variants are ahead of the WoLF variants throughout, most clearly for the chaotic Lorenz-63 model, although adjusting the learning rate of the WoLF variants can likely recover part of this difference. Heuristically, for the filter to remain stable the observations need to convey enough information for the analysis step to counteract the error growth in the forecast. To a degree, this is contained in \cref{thm: stability} as for the regular KF the steady-state analysis covariance balances both quantities, forecast error and analysis correction. For the DSM KF, the steady-state covariance is replaced by a weakly stochastically bounded analysis covariance. When information is already barely sufficient to account for non-linearity in the signal process, discarding information with no way to account for this loss in some other way can be detrimental. Methods based strictly on detect-and-reject schemes are therefore only applicable when reliable observations are plentiful.

The discussed methods based on GBI employ an adjusted information processing. The WoLF KF in \cite{duran2024outlier} dynamically inflates the observation covariance matrix. The DSM KF employs a different discrepancy measure, which for the presented choice of diffusion matrix and weight kernel results in an observation covariance rescaling with factor in $(0,2]$ plus the observation shift. This increases information intake when forecast and observation align well, as indicated by the theoretical results in \cref{sec: stability,sec: tuning}, although rigorous only for the linear case. However, even with adjusted information intake in either direction, i.e., dynamically adjusting the correction in the analysis step, a sufficient amount of information needs to be available likely requiring a stronger assumption than the bounded innovations in \cref{thm: stability} for highly non-linear system dynamics.

\section{Conclusion}

This paper introduced the DSM KF employing diffusion score matching for generalised Bayesian filtering. It combines closed-form Gaussian updates with global bias robustness, a stability guarantee for the analysis covariance, and a default tuning. We derived its EnKF, ESRF and LETKF variants for non-linear dynamical systems. In simulation experiments, the DSM filter variants remained stable with good CRPS under contamination where the regular filters lost track. With this, we aim to further the discussion about generalised Bayesian inference in data assimilation and contribute to a shared language between both communities.

The discussion in the work at hand is in no way exhaustive or complete. Especially for considering non-linear dynamical systems and ensemble approximations, questions on accuracy such as in \cite{takeda2024uniform,takeda2025quantifying} for the regular EnKF are crucial. A central challenge lies in that for the DSM KF variants, additional tools are required for controlling stochasticity and dependence introduced via weight kernels. The result on stability in \cref{thm: stability} and the approach to tuning investigated some of the tools that may be useful in that regard, yet analysis of the analysis mean is only considered at the level of a single update (\cref{lem: influence}). Lyapunov arguments as employed for the analysis of EnKFs \cite{tong2016nonlinear,dewiljes2018long} combined with additional conditions on the true DGPs $\pi$ may provide a first approach to balance forecast error increase and analysis correction for a rigorous error stability analysis. Additionally, non-linear observation operators were not considered in the derivation of the ensemble variants. They provide a natural next question as they may strongly impact the weight kernel.


\begin{funding}
This work has been funded by Deutsche Forschungsgemeinschaft (DFG) - Project-ID 318763901 - SFB1294. Presentation of results at a workshop was partially funded by the European Metrological Society YSCA.
\end{funding}

\begin{disclosure}
    \textbf{Programming.} All code was implemented in \texttt{R} version 4.5.1, with all computationally expensive functions implemented in C++ via \texttt{Rcpp} version 1.1.1 and \texttt{RcppArmadillo} version 15.2.4-1. Repetitions were run in parallel with \texttt{future} version 1.69.0 and \texttt{future.apply} version 1.20.2. The code implementations are elementary and code is available on request.\\
    \textbf{Generative AI.}
    Claude Opus 5.5 by Anthropic was used to revise a previous draft of this work in a final stage. It contributed to the following changes:
    \begin{itemize}
        \item \textbf{Review.} The previous draft was reviewed and a list of flaws including theoretical inconsistencies, ambiguities, oversights was produced 
        \item \textbf{Code.} An earlier code base was partially restructured by translating previous \texttt{R}-code into \texttt{RCPP}, function tests were added and additional evaluation metrics for the supplement were implemented. 
        \item \textbf{Graphics.} Additional components were added to the \texttt{ggplot} functions of earlier graph implementations.
        \item \textbf{Theoretical results.} The theoretical results in \cref{sec: stability,sec: tuning} were adjusted from earlier versions. In \cref{sec: stability}, the results were adjusted to lose an additional assumption required for directly using the result in \cite{solo1996stability}. In \cref{sec: tuning}, the upper bound on the Jensen gap was sharpened in the constant.
        \item \textbf{Language editing.} Spelling and grammatical errors, notation inconsistencies, and ambiguous wording were pointed out as well as wording was revised. 
        \item \textbf{Literature resources.} Additional literature was suggested with some included in the final draft.
    \end{itemize}
    All AI-assisted content was checked by the authors, e.g. code was validated against previous test cases, every change to the research paper was adjusted and implemented by the authors only.
\end{disclosure}


\bibliographystyle{ba2}
\bibliography{refrences}

@mastersthesis{reimann2024towards,
  title={Towards robust inference for {B}ayesian filtering of linear {G}aussian dynamical systems subject to additive change},
  author={Reimann, Hans},
  year={2024},
  school={Universit{\"a}t Potsdam}
}

@book{sarkka2023bayesian,
  title={{B}ayesian filtering and smoothing},
  author={S{\"a}rkk{\"a}, Simo and Svensson, Lennart},
  volume={17},
  year={2023},
  publisher={Cambridge University Press}
}

@book{reich2015probabilistic,
  title={Probabilistic forecasting and {B}ayesian data assimilation},
  author={Reich, Sebastian and Cotter, Colin},
  year={2015},
  publisher={Cambridge University Press}
}

@book{anderson2012optimal,
  title={Optimal filtering},
  author={Anderson, Brian DO and Moore, John B},
  year={2012},
  publisher={Courier Corporation}
}

@article{zellner1988optimal,
  title={Optimal information processing and {B}ayes's theorem},
  author={Zellner, Arnold},
  journal={The American Statistician},
  volume={42},
  number={4},
  pages={278--280},
  year={1988},
  publisher={Taylor \& Francis}
}

@article{provost2023adaptive,
  title={An adaptive ensemble filter for heavy-tailed distributions: tuning-free inflation and localization},
    author = {Le Provost, M. and Baptista, Ricardo and Eldredge, Jeff and Marzouk, Youssef},
  journal={arXiv preprint arXiv:2310.08741},
  year={2023},
  doi = {10.48550/arXiv.2310.08741}
}

@InProceedings{altamirano2023robust,
  title = 	 {Robust and Scalable {B}ayesian Online Changepoint Detection},
  author =       {Altamirano, Matias and Briol, Francois-Xavier and Knoblauch, Jeremias},
  booktitle = 	 {Proceedings of the 40th International Conference on Machine Learning},
  pages = 	 {642--663},
  year = 	 {2023},
  editor = 	 {Krause, Andreas and Brunskill, Emma and Cho, Kyunghyun and Engelhardt, Barbara and Sabato, Sivan and Scarlett, Jonathan},
  volume = 	 {202},
  series = 	 {Proceedings of Machine Learning Research},
  month = 	 {23--29 Jul},
  publisher =    {PMLR},
  url = 	 {https://proceedings.mlr.press/v202/altamirano23a.html}
}

@InProceedings{altamirano2023robust2,
  title = 	 {Robust and Conjugate {G}aussian Process Regression},
  author =       {Altamirano, Matias and Briol, Francois-Xavier and Knoblauch, Jeremias},
  booktitle = 	 {Proceedings of the 41st International Conference on Machine Learning},
  pages = 	 {1155--1185},
  year = 	 {2024},
  editor = 	 {Salakhutdinov, Ruslan and Kolter, Zico and Heller, Katherine and Weller, Adrian and Oliver, Nuria and Scarlett, Jonathan and Berkenkamp, Felix},
  volume = 	 {235},
  series = 	 {Proceedings of Machine Learning Research},
  month = 	 {21--27 Jul},
  publisher =    {PMLR},
  url = 	 {https://proceedings.mlr.press/v235/altamirano24a.html}
}

@article{boustati2020generalised,
  title={Generalised {B}ayesian filtering via sequential {M}onte {C}arlo},
  author={Boustati, Ayman and Akyildiz, Omer Deniz and Damoulas, Theodoros and Johansen, Adam},
  journal={Advances in Neural Information Processing Systems},
  volume={33},
  pages={418--429},
  year={2020}
}

@article{barp2019minimum,
  title={Minimum {S}tein discrepancy estimators},
  author={Barp, Alessandro and Briol, Fran{\c{c}}ois-Xavier and Duncan, Andrew and Girolami, Mark and Mackey, Lester},
  journal={Advances in Neural Information Processing Systems},
  volume={32},
  year={2019}
}

@article{ferro2014fair,
  title={Fair scores for ensemble forecasts},
  author={Ferro, Christopher AT},
  journal={Quarterly Journal of the Royal Meteorological Society},
  volume={140},
  number={683},
  pages={1917--1923},
  year={2014},
  publisher={Wiley Online Library}
}

@book{evensen2009data,
  title={Data assimilation: the ensemble Kalman filter},
  author={Evensen, Geir},
  year={2009},
  publisher={Springer}
}

@article{jewson2018principles,
  title={Principles of {B}ayesian inference using general divergence criteria},
  author={Jewson, Jack and Smith, Jim Q and Holmes, Chris},
  journal={Entropy},
  volume={20},
  number={6},
  pages={442},
  year={2018},
  publisher={MDPI}
}

@article{bissiri2016general,
  title={A general framework for updating belief distributions},
  author={Bissiri, Pier Giovanni and Holmes, Chris C and Walker, Stephen G},
  journal={Journal of the Royal Statistical Society Series B: Statistical Methodology},
  volume={78},
  number={5},
  pages={1103--1130},
  year={2016},
  publisher={Oxford University Press}
}

@article{matsubara2022robust,
  title={Robust generalised {B}ayesian inference for intractable likelihoods},
  author={Matsubara, Takuo and Knoblauch, Jeremias and Briol, Fran{\c{c}}ois-Xavier and Oates, Chris J},
  journal={Journal of the Royal Statistical Society Series B: Statistical Methodology},
  volume={84},
  number={3},
  pages={997--1022},
  year={2022},
  publisher={Oxford University Press}
}

@article{hyvarinen2005estimation,
  author  = {Aapo Hyv{{\"a}}rinen},
  title   = {Estimation of Non-Normalized Statistical Models by Score Matching},
  journal = {Journal of Machine Learning Research},
  year    = {2005},
  volume  = {6},
  number  = {24},
  pages   = {695--709},
  url     = {http://jmlr.org/papers/v6/hyvarinen05a.html}
}

@book{huber2004robust,
  title={Robust statistics},
  author={Huber, Peter J},
  volume={523},
  year={2004},
publisher={John Wiley \& Sons}
}

@article{ghosh2016robust,
  title={Robust {B}ayes estimation using the density power divergence},
  author={Ghosh, Abhik and Basu, Ayanendranath},
  journal={Annals of the Institute of Statistical Mathematics},
  volume={68},
  pages={413--437},
  year={2016},
  publisher={Springer}
}

@book{evensen2022data,
  title={Data assimilation fundamentals: A unified formulation of the state and parameter estimation problem},
  author={Evensen, Geir and Vossepoel, Femke C and van Leeuwen, Peter Jan},
  year={2022},
  publisher={Springer Nature}
}

@article{tang2024generalized,
  title={A Generalized t-Distribution-Based Kernel Adaptive Filtering Algorithm},
  author={Tang, Huchuan and Han, Hongyu and Zhang, Sheng and Feng, Wenting},
  journal={IEEE Transactions on Circuits and Systems II: Express Briefs},
  year={2024},
  volume={71},
  number={6},
  pages={3241-3245},
  publisher={IEEE},
  doi={10.1109/TCSII.2024.3356912}}

@article{mahalanobis2018generalized,
	title = {On the generalized distance in statistics},
	journal = {Proceedings of the National Institute of Science of India},
	author = {Mahalanobis, P.C.},
	year = {1936},
    volume={2},
	pages = {49--55},
    publisher={National Institute of Science of India}
}

@inproceedings{solo1996stability,
  title={Stability of the {K}alman filter with stochastic time-varying parameters},
  author={Solo, Victor},
  booktitle={Proceedings of 35th IEEE Conference on Decision and Control},
  volume={1},
  pages={57--61},
  year={1996},
  organization={IEEE}
}

@inproceedings{knoblauch2018spatio,
  title={Spatio-temporal {B}ayesian on-line changepoint detection with model selection},
  author={Knoblauch, Jeremias and Damoulas, Theodoros},
  booktitle={International Conference on Machine Learning},
  pages={2718--2727},
  year={2018},
  organization={PMLR}
}

@InProceedings{duran2024outlier,
  title = 	 {Outlier-robust Kalman Filtering through Generalised {B}ayes},
  author =       {Duran-Martin, Gerardo and Altamirano, Matias and Shestopaloff, Alex and S\'{a}nchez-Betancourt, Leandro and Knoblauch, Jeremias and Jones, Matt and Briol, Francois-Xavier and Murphy, Kevin Patrick},
  booktitle = 	 {Proceedings of the 41st International Conference on Machine Learning},
  pages = 	 {12138--12171},
  year = 	 {2024},
  editor = 	 {Salakhutdinov, Ruslan and Kolter, Zico and Heller, Katherine and Weller, Adrian and Oliver, Nuria and Scarlett, Jonathan and Berkenkamp, Felix},
  volume = 	 {235},
  series = 	 {Proceedings of Machine Learning Research},
  month = 	 {21--27 Jul},
  publisher =    {PMLR},
  url = 	 {https://proceedings.mlr.press/v235/duran-martin24a.html}
}

@article{hunt2007efficient,
  title={Efficient data assimilation for spatiotemporal chaos: A local ensemble transform {K}alman filter},
  author={Hunt, Brian R and Kostelich, Eric J and Szunyogh, Istvan},
  journal={Physica D: Nonlinear Phenomena},
  volume={230},
  number={1-2},
  pages={112--126},
  year={2007},
  publisher={Elsevier}
}

@book{kullback1997information,
  title={Information Theory and Statistics},
  author={Kullback, S.},
  lccn={59005881},
  series={Wiley publication in mathematical statistics},
  url={https://books.google.de/books?id=XeRQAAAAMAAJ},
  year={1959},
  publisher={Wiley}
}

@article{gong2021interpreting,
  title={Interpreting diffusion score matching using normalizing flow},
  author={Gong, Wenbo and Li, Yingzhen},
  journal={arXiv preprint arXiv:2107.10072},
  year={2021}
}

@inproceedings{mcallester1998some,
  title={Some {PAC}-{B}ayesian theorems},
  author={McAllester, David A},
  booktitle={Proceedings of the eleventh annual Conference on Computational Learning Theory},
  pages={230--234},
  year={1998}
}

@article{takeda2024uniform,
  title={Uniform Error Bounds of the Ensemble Transform {K}alman Filter for Chaotic Dynamics with Multiplicative Covariance Inflation},
  author={Takeda, Kota and Sakajo, Takashi},
  journal={SIAM/ASA Journal on Uncertainty Quantification},
  volume={12},
  number={4},
  pages={1315--1335},
  year={2024},
  publisher={SIAM}
}

@article{gao2017bounds,
  title={Bounds on the {J}ensen gap, and implications for mean-concentrated distributions},
  author={Gao, Xiang and Sitharam, Meera and Roitberg, Adrian E},
  journal={arXiv preprint arXiv:1712.05267},
  year={2017}
}

@article{yatawara1991kalman,
  title={A {K}alman filter in the presence of outliers},
  author={Yatawara, Nihal and Abraham, Bovas and MacGregor, John F},
  journal={Communications in Statistics-Theory and Methods},
  volume={20},
  number={5-6},
  pages={1803--1820},
  year={1991},
  publisher={Taylor \& Francis}
}

@article{xie1994robust,
  title={Robust {K}alman filtering for uncertain systems},
  author={Xie, Lihua and Soh, Yeng Chai},
  journal={Systems \& Control Letters},
  volume={22},
  number={2},
  pages={123--129},
  year={1994},
  publisher={Elsevier}
}

@article{sun2022control,
  title={Control simulation experiments of extreme events with the Lorenz-96 model},
  author={Sun, Qiwen and Miyoshi, Takemasa and Richard, Serge},
  journal={Nonlinear Processes in Geophysics},
  volume={30},
  number={2},
  pages={117--128},
  year={2023},
  publisher={Copernicus Publications G{\"o}ttingen, Germany}
}

@article{takeda2025quantifying,
  title={Noise-scaled accuracy of the ensemble Kalman filter with an instability-based minimum ensemble size},
  author={Takeda, Kota and Miyoshi, Takemasa},
  journal={Nonlinear Processes in Geophysics},
  volume={33},
  number={3},
  pages={335--346},
  year={2026},
  publisher={Copernicus Publications G{\"o}ttingen, Germany}
}

@article{mu2023huber,
  title={On {H}uber's contaminated model},
  author={Mu, Weiyan and Xiong, Shifeng},
  journal={Journal of Complexity},
  volume={77},
  pages={101745},
  year={2023},
  publisher={Elsevier}
}

@article{tong2016nonlinear,
  title={Nonlinear stability and ergodicity of ensemble based Kalman filters},
  author={Tong, Xin T and Majda, Andrew J and Kelly, David},
  journal={Nonlinearity},
  volume={29},
  number={2},
  pages={657--691},
  year={2016},
  publisher={IOP Publishing}
}

@article{dewiljes2018long,
  title={Long-time stability and accuracy of the ensemble Kalman--Bucy filter for fully observed processes and small measurement noise},
  author={De Wiljes, Jana and Reich, Sebastian and Stannat, Wilhelm},
  journal={SIAM Journal on Applied Dynamical Systems},
  volume={17},
  number={2},
  pages={1152--1181},
  year={2018},
  publisher={SIAM}
}

@article{gneiting2007strictly,
  title={Strictly proper scoring rules, prediction, and estimation},
  author={Gneiting, Tilmann and Raftery, Adrian E},
  journal={Journal of the American Statistical Association},
  volume={102},
  number={477},
  pages={359--378},
  year={2007},
  publisher={Taylor \& Francis}
}

@article{supplement,
  title={Supplement to "Robust Ensemble Kalman filtering Under Observation Noise Misspecification"},
  author={Reimann, Hans and Reich, Sebastian},
  year={2026},
  doi = {TBD}
}

@article{fortin2014should,
  title={Why should ensemble spread match the RMSE of the ensemble mean?},
  author={Fortin, Vincent and Abaza, Mabrouk and Anctil, Francois and Turcotte, Raphael},
  journal={Journal of Hydrometeorology},
  volume={15},
  number={4},
  pages={1708--1713},
  year={2014}
}

@article{evensen1994sequential,
  title={Sequential data assimilation with a nonlinear quasi-geostrophic model using Monte Carlo methods to forecast error statistics},
  author={Evensen, Geir},
  journal={Journal of Geophysical Research: Oceans},
  volume={99},
  number={C5},
  pages={10143--10162},
  year={1994},
  publisher={Wiley Online Library}
}

@article{houtekamer2016review,
  title={Review of the ensemble Kalman filter for atmospheric data assimilation},
  author={Houtekamer, Peter L and Zhang, Fuqing},
  journal={Monthly Weather Review},
  volume={144},
  number={12},
  pages={4489--4532},
  year={2016}
}

@article{van2019particle,
  title={Particle filters for high-dimensional geoscience applications: A review},
  author={Van Leeuwen, Peter Jan and K{\"u}nsch, Hans R and Nerger, Lars and Potthast, Roland and Reich, Sebastian},
  journal={Quarterly Journal of the Royal Meteorological Society},
  volume={145},
  number={723},
  pages={2335--2365},
  year={2019},
  publisher={Wiley Online Library}
}

\clearpage

\setcounter{section}{0}
\setcounter{figure}{0}
\renewcommand{\thesection}{S\arabic{section}}
\renewcommand{\thefigure}{S\arabic{figure}}

\makeatletter
\begingroup
  \frontmatter@style                   
  \vglue\title@skip                    
  \nointerlineskip                     
  {\title@style\title@size             
   Supplementary Material for ``Robust Ensemble Kalman Filtering Under Observation Noise Misspecification''\par}
  \vskip\authors@skip                  
  {\noindent\author@style\author@size  
   Hans Reimann and Sebastian Reich\par}
\endgroup
\frontmatter@cmd                       
\makeatother

\begingroup
\begin{frontmatter}
\supptitle{}
\runtitle{Supplement: Robust EnKF Under Observation Misspecification}
\end{frontmatter}
\endgroup

This supplement provides example trajectories with Gaussian approximation credible bands as well as additional evaluation metrics for all simulation experiments in \cref{sec: simulation study}, i.e. the time-averaged RMSE, the log-score, the empirical coverage of the $90\%$ credible intervals, and the spread-RMSE ratio. Additionally, it provides details on the implementation of inflation and localisation for the Lorenz-96 simulation experiment in \cref{sec: simulation study NL}.
\cref{sec: simulation study LGSS}

\section{Evaluation metrics}
In addition to the CRPS of \cref{sec: simulation study}, we report additional evaluation metrics. Let $x^\dagger$ denote the true state, $m^a$ denote the analysis mean, and $P^a$ denote the analysis covariance matrix. For the EnKF and LETKF the analysis mean and covariance are replaced by their empirical counterparts. Evaluation of the LETKF just considers the diagonal entries of the empirical covariance matrix. The RMSE, log-score, and empirical coverage are averaged over the analysis time window. For all metrics the median over the $M_\mathrm{MC}$ repeated runs is reported. Runs losing track of the signal are included in the evaluation. We consider the following:
\begin{itemize}
    \item The RMSE (lower is better),
    \begin{equation*}
        \mathrm{RMSE}(x^\dagger;\, m^a)\coloneqq\sqrt{\frac{1}{d_X}\sum_{i=1}^{d_X}\left(x^\dagger_i-m^a_i\right)^2}\,,
    \end{equation*}
    \item the log-score \cite{gneiting2007strictly}, i.e. the negative log-density of the Gaussian posterior approximation at the true state, (lower is better),
    \begin{equation*}
        \log\mathrm{S}(x^\dagger;\, m^a,P^a)\coloneqq-\log n(x^\dagger; m^a, P^a)\,,
    \end{equation*}
    \item the empirical coverage of the componentwise 90\% credible intervals (closer to $0.9$ is better, $<0.9$ is overconfident),
    \begin{equation*}
        \mathrm{C}_{90}(x^\dagger;\, m^a,P^a)\coloneqq\frac{1}{d_X}\sum_{i=1}^{d_X}\mathbf{1}\left\{\left|x^\dagger_i-m^a_i\right|\leq \sqrt{P^a_{ii}}z_{0.95}\right\}
    \end{equation*}
    with $z_{0.95}=\Phi^{-1}(0.95)$ and $\Phi$ the standard Gaussian CDF, and
    \item the spread-RMSE ratio \cite{fortin2014should}, i.e. a method's estimate of its RMSE compared to true RMSE, (closer to $1$ is better, $<1$ is overconfident),
    \begin{equation*}
        \mathrm{SER}\coloneqq\sqrt{\frac{\sum_{n}\mathrm{Sp}(P^a_n)^2}{\sum_{n}\mathrm{RMSE}(x^\dagger_n;\, m^a_n)^2}}
        \quad\text{with}\quad
        \mathrm{Sp}(P^a)\coloneqq\sqrt{\frac{1}{d_X}\sum_{i=1}^{d_X} P^a_{ii}}\,.
    \end{equation*}
\end{itemize}

\section{Additional simulation results}

\subsection{Target tracking}\label{sec: S TT}

\begin{figure}[h]
  \centering
  \includegraphics[width=\linewidth]{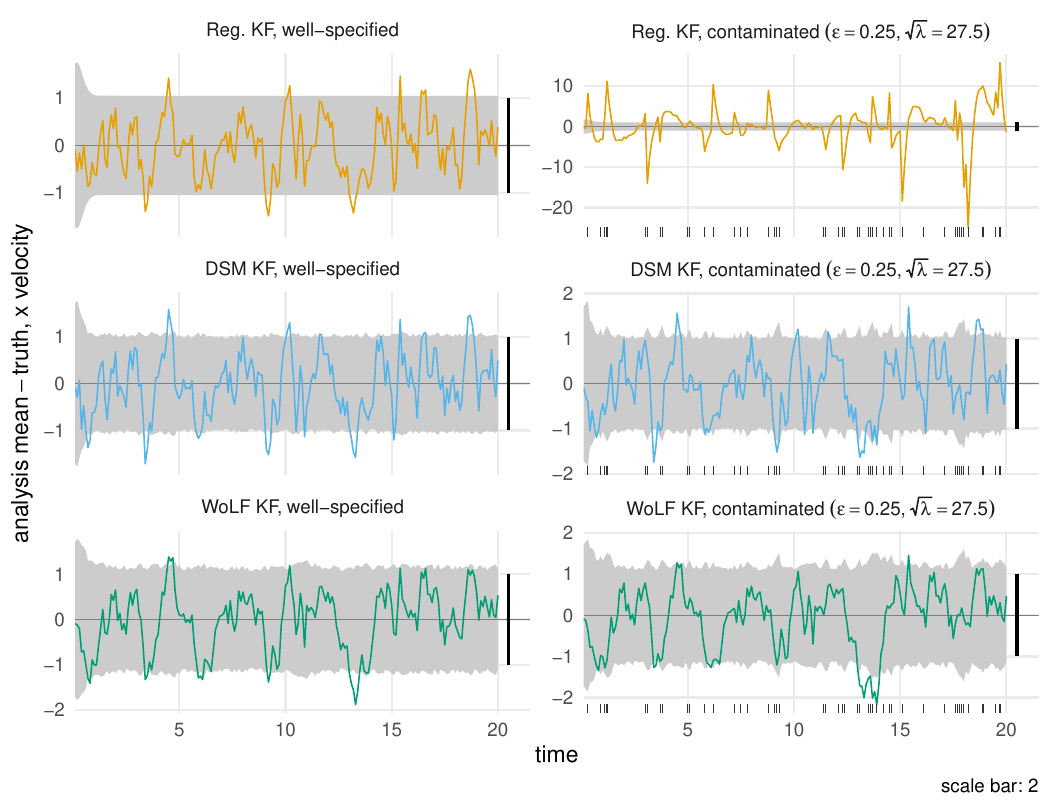}
  \caption{Error of the analysis mean in the unobserved first velocity component with
    the 90\% credible band (grey) of the regular, DSM and WoLF KF (top to bottom) for one
    repetition, well specified (left) and contaminated with $\varepsilon = 0.25$,
    $\sqrt{\lambda} = 27.5$ (right), over the first 200 evaluated analysis times. Ticks mark
    contaminated observations. The black bar has the same length in all panels.}
  \label{fig: S traj TT}
\end{figure}
 
\begin{figure}[p]
  \centering
  \includegraphics[width=\linewidth]{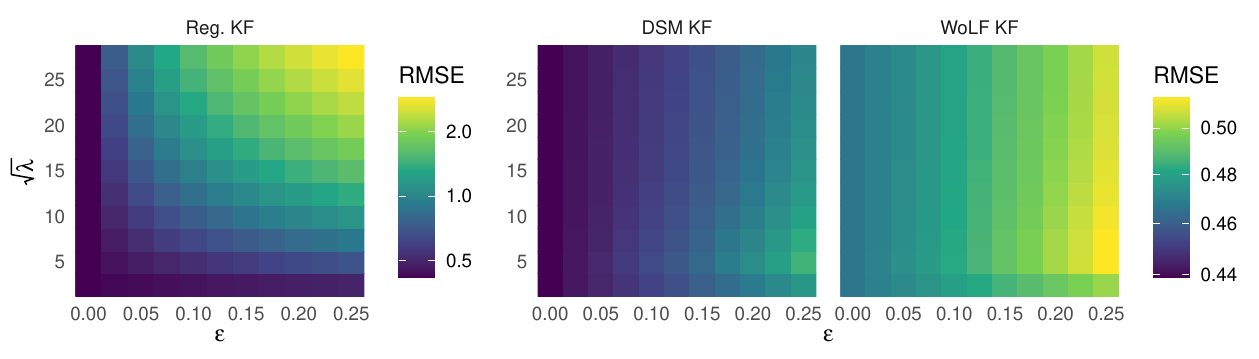}\\[4pt]
  \includegraphics[width=\linewidth]{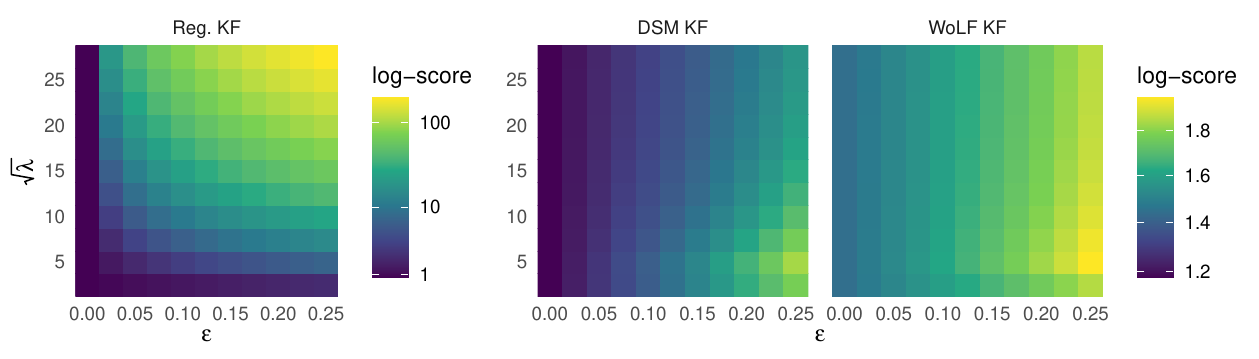}\\[4pt]
  \includegraphics[width=\linewidth]{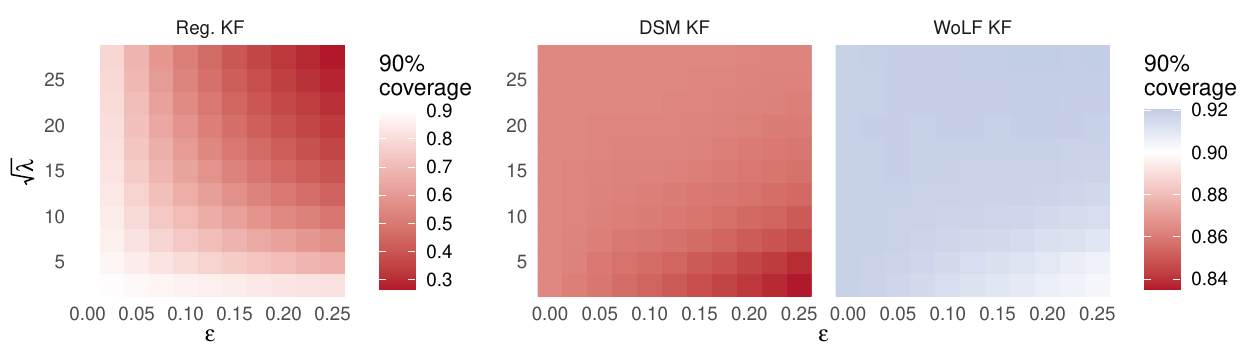}\\[4pt]
  \includegraphics[width=\linewidth]{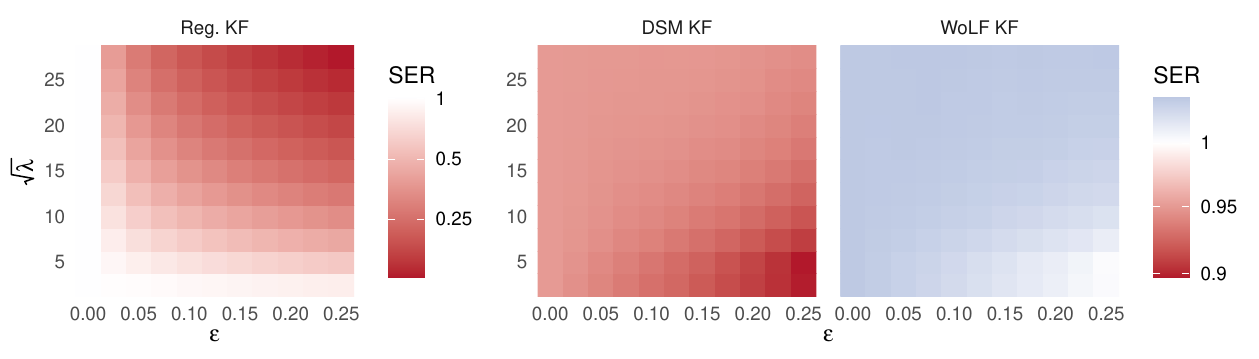}
  \caption{2D target tracking (\cref{sec: simulation study LGSS}): median RMSE, log-score,
    90\% coverage, and spread-RMSE ratio (top to bottom) over $M_{\mathrm{MC}} = 2500$ repetitions of the regular KF
    (left), DSM KF (middle) and WoLF KF (right) for different frequencies ($\varepsilon$) and
    degrees ($\lambda$). The DSM and WoLF KF share the right scale. Coverage and spread-RMSE ratio are white at their respective ideal values $0.9$ and $1$.}
  \label{fig: S heat TT}
\end{figure}

\clearpage
\subsection{Stochastic Lorenz-63 model}\label{sec: S L63}

\begin{figure}[h]
  \centering
  \includegraphics[width=\linewidth]{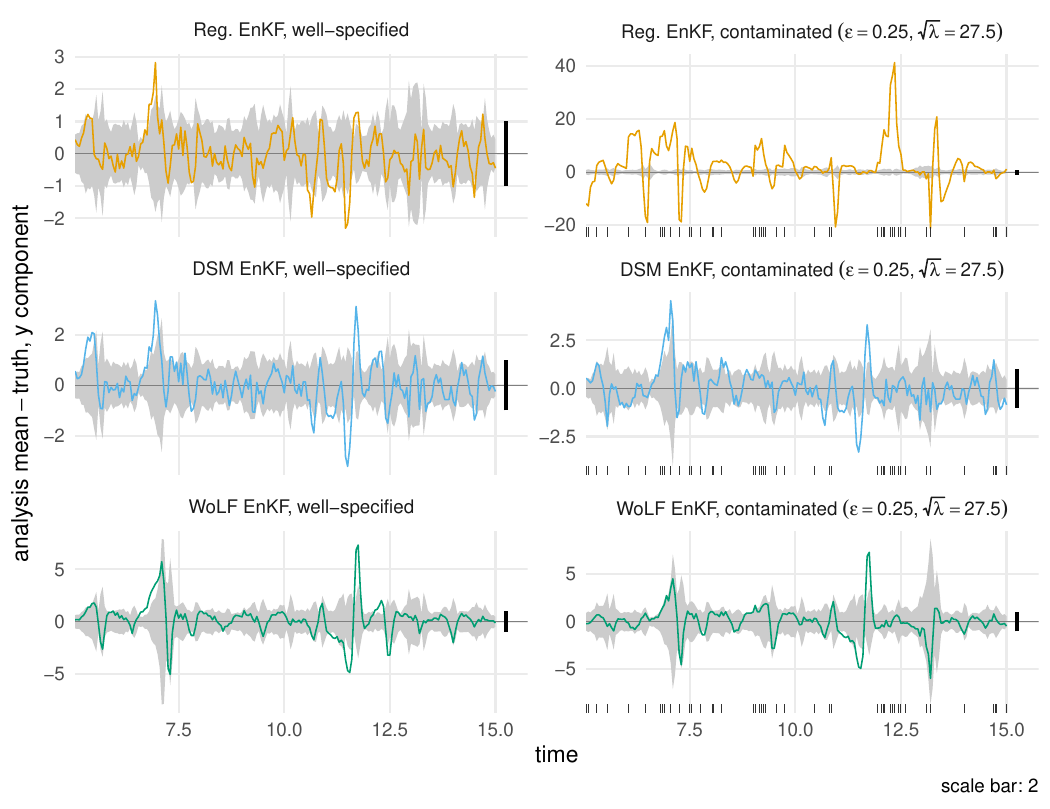}
  \caption{Error of the analysis mean in the unobserved second state component with
    the 90\% credible band (grey) of the regular, DSM and WoLF EnKF (top to bottom) for one
    repetition, well specified (left) and contaminated with $\varepsilon = 0.25$,
    $\sqrt{\lambda} = 27.5$ (right), over the first 200 evaluated analysis times. Ticks mark
    contaminated observations. The black bar has the same length in all panels.}
  \label{fig: S traj L63}
\end{figure}

\begin{figure}[p]
  \centering
  \includegraphics[width=\linewidth]{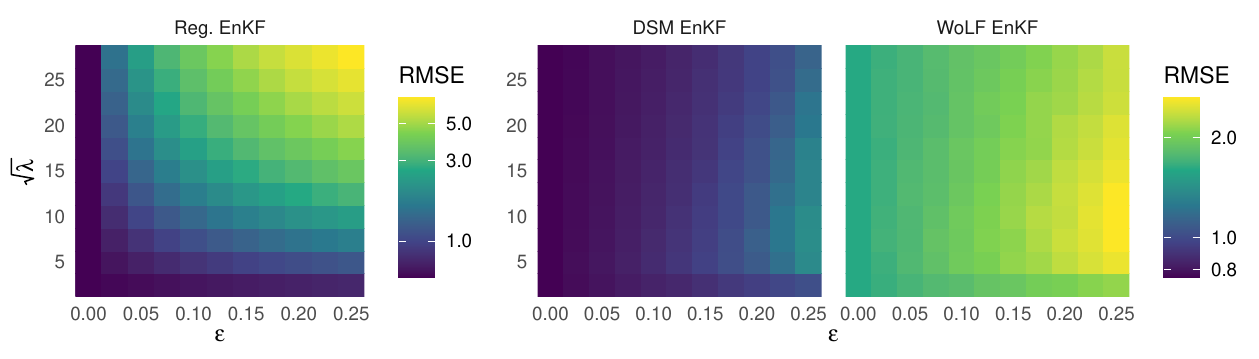}\\[4pt]
  \includegraphics[width=\linewidth]{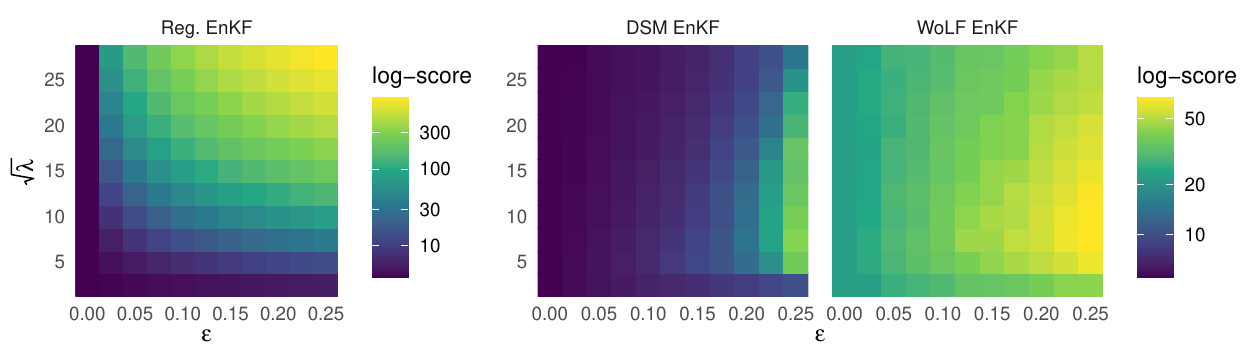}\\[4pt]
  \includegraphics[width=\linewidth]{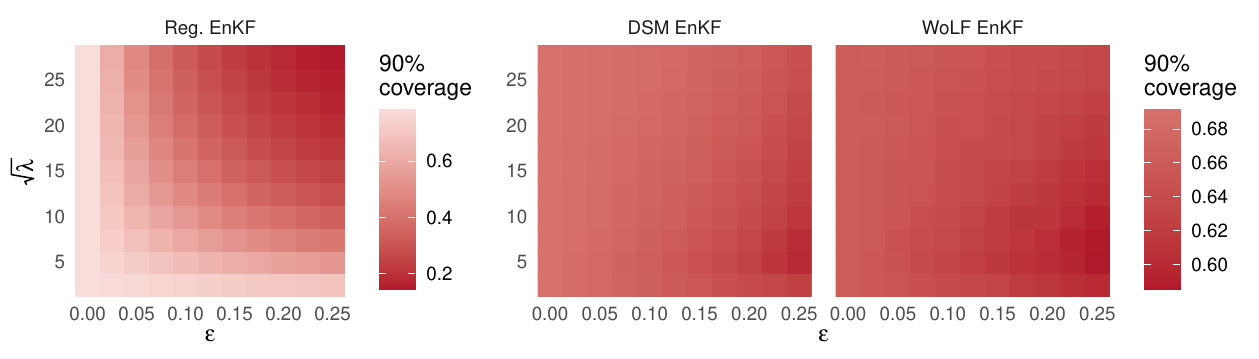}\\[4pt]
  \includegraphics[width=\linewidth]{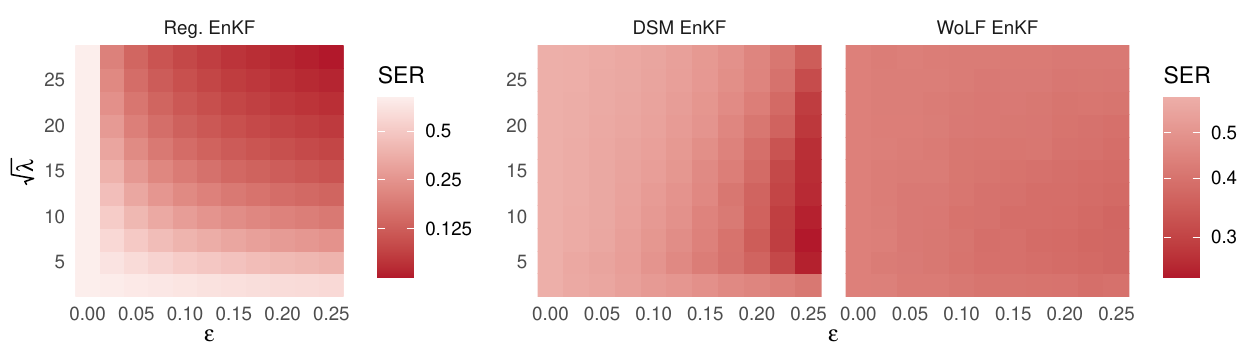}
  \caption{Stochastic Lorenz-63 model (\cref{sec: simulation study NL}): median RMSE, log-score,
    90\% coverage, and spread-RMSE ratio (top to bottom) over $M_{\mathrm{MC}} = 1000$ repetitions of the regular EnKF
    (left), DSM EnKF (middle) and WoLF EnKF (right) for different frequencies ($\varepsilon$) and
    degrees ($\lambda$). The DSM and WoLF EnKF share the right scale. Coverage and spread-RMSE ratio are white at their respective ideal values $0.9$ and $1$.}
  \label{fig: S heat L63}
\end{figure}

\clearpage
\subsection{40D Lorenz-96 model}\label{sec: S L96}

Covariance inflation and localisation are implemented as stated in \cite{sun2022control}. The regular, DSM and WoLF LETKF apply multiplicative covariance inflation to the forecast covariance $P^f$ with a factor $\rho_{\mathrm{inf}}=1.06$. In ensemble subspace this replaces $M-1$ by $\frac{M-1}{\rho_{\mathrm{inf}}}$. Localisation refers to $R$-localisation as introduced in \cite{hunt2007efficient}. The $j$-th component is updated via its own local analysis on the $39$ local points $i=j-19,\ldots,j+19$ (modulo $40$ on the torus) and only component $j$ is kept. The local $R_{\mathrm{loc}}$ is diagonal with entries $\exp\left(\frac{\mathrm{d}(i,j)^2}{2L^2}\right)$ for $\mathrm{d}(i,j)=\min\{|i-j\mod40|,\ |j-i\mod40|\}$ for fixed $j$ and $i=j-19,\ldots,j+19$ and $L=5.45$. In \cite{sun2022control} localisation is tuned to minimise RMSE achieving a value of $0.1989$. The regular LETKF in the well-specified case in our implementation matches this value achieving a median RMSE of $0.198$.

The weight kernels of the DSM and WoLF LETKF are evaluated for each local analysis with the inflated forecast covariance $\rho_{\mathrm{inf}}P^f$ and local observation covariance $R_{\mathrm{loc}}$. The tuning parameter $q^2$ (and $c^2$ in the WoLF LETKF) employs the effective number of observations $d_Y^{\mathrm{eff}}=\sum_{i=j-19}^{j+19}\exp\left(-\frac{\mathrm{d}(i,j)^2}{2L^2}\right)\approx 13.66$ since we assumed $R=\mathbf{I}_{40\times40}$.

\begin{figure}[h]
  \centering
  \includegraphics[width=\linewidth]{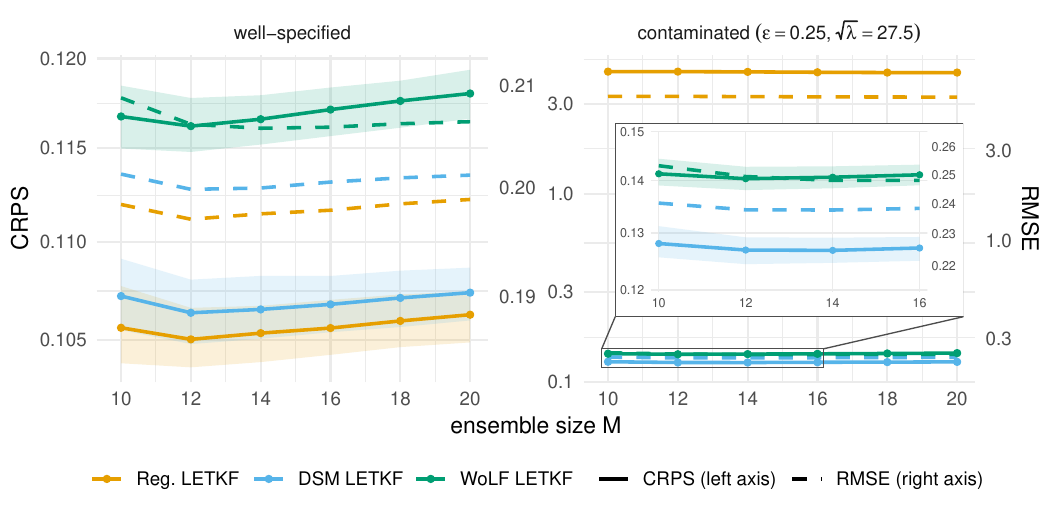}
  \caption{Median CRPS (solid, shaded $25$--$75\%$ quantile band; left $y$-axis) and RMSE (dashed; right $y$-axis) over $M_\mathrm{MC}=1000$ repetitions of LETKF variants for different ensemble sizes in the well-specified (left) and contaminated model (right).}
  \label{fig: S ens L96}
\end{figure}

\begin{figure}[p]
  \centering
  \includegraphics[width=\linewidth]{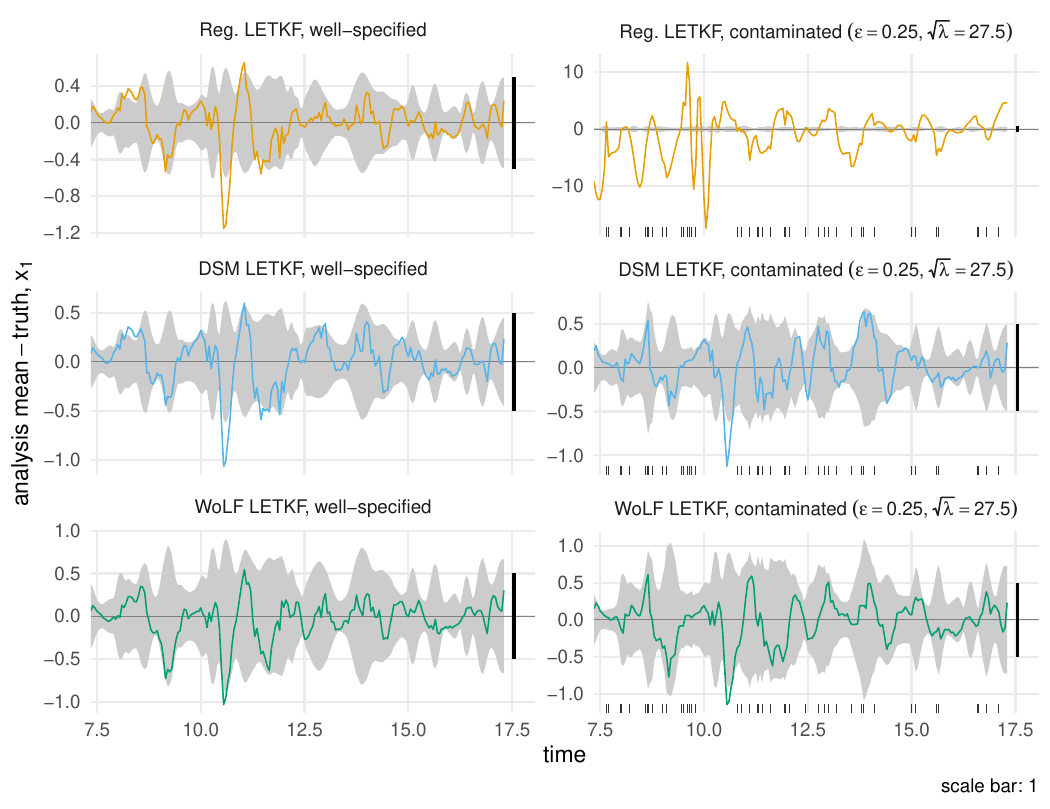}
  \caption{Error of the analysis mean in the first state component with
    the 90\% credible band (grey) of the regular, DSM and WoLF LETKF (top to bottom) for one
    repetition, well specified (left) and contaminated with $\varepsilon = 0.25$,
    $\sqrt{\lambda} = 27.5$ (right), over the first 200 evaluated analysis times. Ticks mark
    contaminated observations. The black bar has the same length in all panels.}
  \label{fig: S traj L96}
\end{figure}

\begin{figure}[p]
  \centering
  \includegraphics[width=\linewidth]{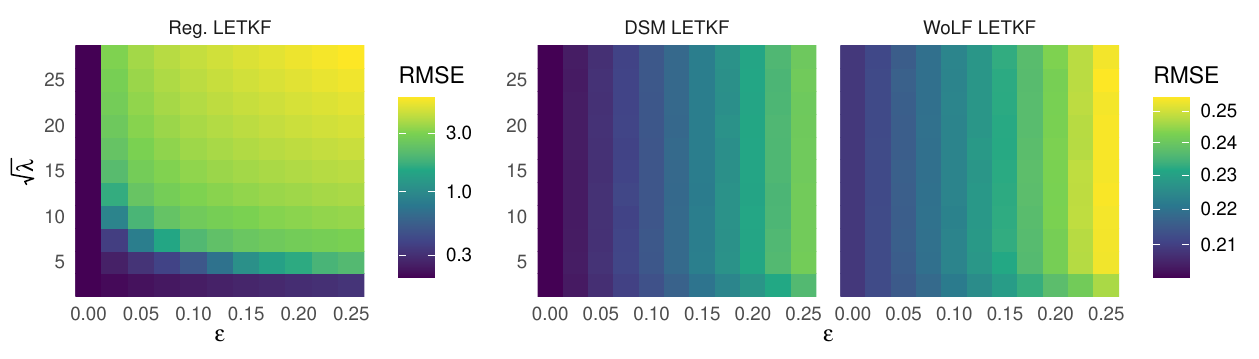}\\[4pt]
  \includegraphics[width=\linewidth]{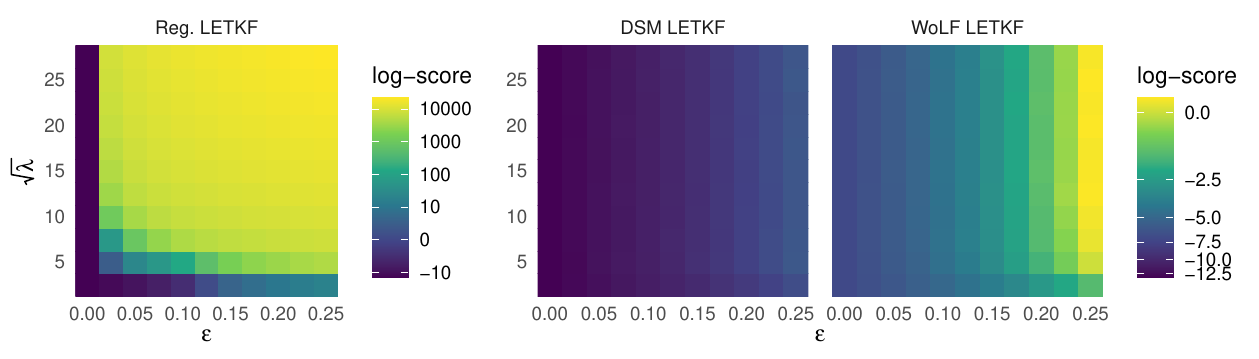}\\[4pt]
  \includegraphics[width=\linewidth]{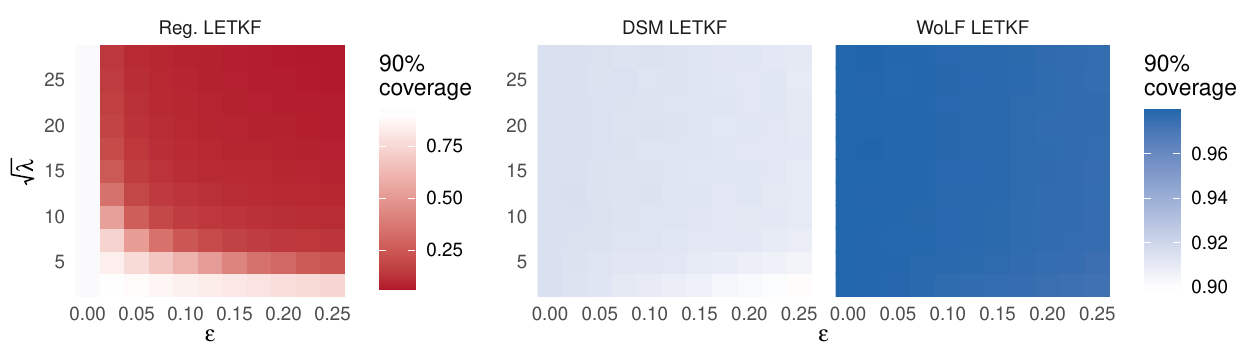}\\[4pt]
  \includegraphics[width=\linewidth]{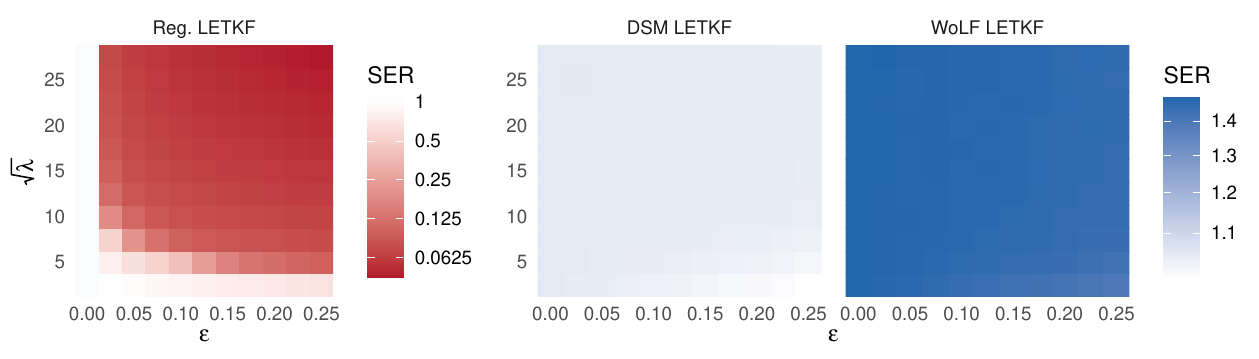}
  \caption{40D Lorenz-96 model (\cref{sec: simulation study NL}): median RMSE, log-score,
    90\% coverage, and spread-RMSE ratio (top to bottom) over $M_{\mathrm{MC}} = 500$ repetitions of the regular LETKF
    (left), DSM LETKF (middle) and WoLF LETKF (right) for different frequencies ($\varepsilon$) and
    degrees ($\lambda$). The DSM and WoLF LETKF share the right scale. Coverage and spread-RMSE ratio are white at their respective ideal values $0.9$ and $1$.}
  \label{fig: S heat L96}
\end{figure}

\end{document}